\documentclass[twoside]{article}

\usepackage[preprint]{aistats2026}
\usepackage[round]{natbib}

\usepackage{booktabs}       
\usepackage{amsfonts}
\usepackage{amsmath}
\usepackage{amsthm}
\usepackage{amssymb}
\usepackage{nicefrac}       
\usepackage{microtype}   

\usepackage{xcolor}         
\usepackage{graphicx}

\usepackage{subcaption}
\usepackage{caption}
\usepackage{comment}

\newtheorem{theorem}{Theorem}

\newtheorem{lemma}{Lemma}
\newtheorem{definition}{Definition}
\newtheorem{corollary}{Corollary}
\newtheorem{remark}{Remark}

\begin{document}

%

%

\twocolumn[

\aistatstitle{Asymptotic Optimality Theory of Confidence Intervals of the Mean}

\aistatsauthor{  Vikas Deep \And Achal Bassamboo \And  Sandeep Juneja }

\aistatsaddress{ NUS, Singapore \And  Kellogg, Northwestern University \And Ashoka University, India} ]

\begin{abstract}
We address the classical problem of constructing confidence intervals (CIs) for the mean of a distribution, given $N$ i.i.d. samples, such that the CI contains the true mean with probability at least $1 - \delta$, where $\delta \in (0,1)$. We characterize three distinct learning regimes based on the minimum achievable limiting width of any CI as the sample size $N_\delta \to \infty$ and $\delta \to 0$. In the first regime, where $N_\delta$ grows slower than $\log(1/\delta)$, the limiting width of any CI equals the width of the distribution’s support, precluding meaningful inference. In the second regime, where $N_\delta$ scales as $\log(1/\delta)$, we precisely characterize the minimum limiting width, which depends on the scaling constant. In the third regime, where $N_\delta$ grows faster than $\log(1/\delta)$, complete learning is achievable, and the limiting width of the CI collapses to zero and CI converges to the true mean. We demonstrate that CIs derived from concentration inequalities based on Kullback-Leibler (KL) divergences achieve asymptotically optimal performance, attaining the minimum limiting width in both the sufficient and the complete learning regimes for distributions in three families: single-parameter exponential, bounded support and known bound on $(1+\epsilon)^{\rm th}$ moment. Additionally, these results extend to one-sided CIs, with the width notion adjusted appropriately. Finally, we generalize our findings to settings with random per-sample costs, motivated by practical applications such as stochastic simulators and cloud service selection. Instead of a fixed sample size, we consider a cost budget $C_\delta$, identifying analogous learning regimes and characterizing the optimal CI construction policy.
\end{abstract}

\section{INTRODUCTION AND MAIN CONTRIBUTION} \label{sec:intro}
The problem of constructing a confidence interval (CI) for the mean of a distribution with a coverage guarantee, ensuring that the mean lies within the CI with probability at least $1-\delta$ for a pre-specified $\delta \in (0,1)$, is well-studied in the statistics literature. This problem has significant applications in A/B testing, experimentation, data analytics, and simulation. Typically, this is achieved using concentration inequalities for the mean, given a sample size of $N$ i.i.d. observations from a probability distribution (see \cite{hoeffding1994probability}, \cite{waudby2024estimating}, \cite{maurer2009empirical}, \cite{audibert2009exploration}, \cite{boucheron2009bernstein}, \cite{catoni2012challenging}, \cite{chen2021generalized} and \cite{bennett1962probability}). In this paper, we primarily aim to address an important element absent from the literature on CI. \textit{While there are various methods to construct a confidence interval for the mean of a distribution, there are no results characterizing the optimal CI with minimum width.} 
In this paper, we characterize three learning regimes based on the minimum limiting width achievable for any CI construction method/policy as $N_\delta \to \infty$ when $\delta \to 0$, under a mild assumption of stability of policies (as defined in Definition \ref{def_stability}). We would like to emphasize that we consider CI construction policies that provide non-asymptotic coverage guarantees for any fixed $N$ and $\delta$. The only aspect where we rely on asymptotics is in defining our notion of optimality. In a nutshell, under this assumption, the limiting width of CI of mean becomes a deterministic constant for any CI construction method making analysis tractable. The three regimes are defined as follows:

1) \textbf{No learning regime: } If $\lim_{\delta \to 0}\frac{N_\delta}{\log(1/\delta)} \to 0$, the limiting width of the CI is the length of the support of the mean for any stable CI construction policies. This implies that no learning is possible as the sample size is not sufficient.

2) \textbf{Sufficient learning regime: } If $\lim_{\delta \to 0}\frac{N_\delta}{\log(1/\delta)} \to k$ for $k \in (0,\infty)$, we have a sharp characterization of the lower bound on limiting width that involves terms with  KL divergences for any stable CI construction policies. This lower bound on limiting width shrinks as $k$ increases. Further, we show that the method $\pi_1$ (see section \ref{sec_kl}) that constructs CI via inverting concentration inequality based on KL divergences matches the lower bound of the limiting width. Hence, this lower bound result on the limiting width is tight as $\delta \to 0$ and $N_\delta \to \infty $. 

3) \textbf{Complete learning regime:} If $\lim_{\delta \to 0}\frac{N_\delta}{\log(1/\delta)} \to \infty$, the limiting CI width trivially converges to zero. Furthermore, we prove that the method $\pi_1$ has a zero limiting width. In this regime, we further analyze the rate at which the width converges to zero and demonstrate that $\pi_1$ achieves the fastest rate of convergence under some technical assumptions.

Our results apply to both parametric and non-parametric families of probability distributions. For parametric cases, we consider a canonical single-parameter exponential family, while for non-parametric cases, we consider two such settings. The first one is the family of probability distributions with bounded support and the second one is the family of probability distributions with $(1+\epsilon)^{\rm th}$ moment bounded and the bound is known. We assume only that the unknown underlying distribution belongs to this known family. If a CI construction method, that matches the lower bound on the limiting width in the sufficient learning regime and the complete learning regime, is called an asymptotically optimal CI construction method. Next, we state the recipe for asymptotically optimal CI construction policies.

\subsection{Recipe for asymptotically optimal CI construction methods}


To construct an asymptotically optimal CI construction method, when the probability distribution belongs to a canonical single-parametric exponential family, we utilize existing concentration bound in the literature on the deviation of the maximum likelihood estimator (MLE) of the mean from the true mean, measured in terms of the KL divergence (see Eq. \eqref{CI_concen}). Inverting this concentration inequality to construct a CI for the mean leads to an asymptotically optimal CI construction method. For the non-parametric setting, \cite{agrawal2022bandits} and \cite{orabona2023tight} provide a confidence interval that  is constructed via inverting a concentration inequality based on KL divergences. In this setting, we prove the optimality of the CI provided by them. In summary, we prove that the CIs constructed via inverting concentration inequalities based on KL divergences are the best CIs in a reasonable asymptotic regime. 
\subsection{Random sampling cost }
The above results assume that the cost of each sample is fixed and equal to one. However, in many practical situations, this assumption may not hold. For instance, in simulation studies, the cost can represent the time required to simulate the performance of a system design, which may vary depending on system complexity and randomness (see \cite{chick2001new}). Similarly, in real-world applications, online bidding optimization in sponsored search (see \cite{amin2012budget}, \cite{xia2015thompson} and \cite{tran2014efficient}), the cost of acquiring each sample can fluctuate due to factors such as network delays, resource availability, or dynamic market conditions. To include such cases, we extend our results to a setting where each sample is associated with a random cost drawn from an i.i.d. cost distribution. In this scenario, instead of having a fixed sample size $N_\delta$, the constraint is a cost budget $C_\delta$. We similarly identify three corresponding learning regimes and characterize the asymptotically optimal method for this setting as well. When both the rewards and costs of obtaining samples are random, the large deviations of the observed average reward for a given cost budget typically depend on the distributions of both variables (see \cite{glynn2013asymptotic}). However, in our analysis, we find that the limiting width of the CI depends only on the mean of the cost distribution and is otherwise invariant to the cost distribution. 

In summary, our main contributions are as follows:
\begin{itemize}
    \item We characterize three distinct learning regimes: no learning, sufficient learning, and complete learning, based on the relative scaling of the sample size $N_\delta$ with respect to the desired accuracy $1 - \delta$. These regimes lead to markedly different minimum asymptotic limiting widths for any CI construction method, under a mild policy stability assumption.

    \item We derive sharp lower bounds on the limiting CI width in the sufficient learning regime and demonstrate that CIs constructed by inverting existing concentration inequalities involving KL divergences achieve these bounds.This result extends to the complete learning regime, where the CI width converges to zero, establishing that KL divergence-based CI construction is asymptotically optimal in our setting. In the complete learning regime, we further analyze the rate at which the width converges to zero and show that our proposed policy achieves the fastest rate of convergence under some technical assumptions. We emphasize that while CIs constructed by inverting KL-based concentration bounds are known in the literature, their optimality was not established. A key contribution of our work is demonstrating this optimality.
    

    \item We extend our results to more general settings, including random sampling costs, one-sided CIs, and non-parametric distributions. In each of these cases, we again characterize the three learning regimes and an asymptotically optimal CI construction policy. Intriguingly, in the random sampling cost setting, we find that the limiting CI width depends solely on the mean of the cost distribution.
\end{itemize}



\section{LITERATURE REVIEW} \label{sec_lit}
There is a well-established duality between hypothesis testing and the construction of confidence intervals/regions (see \cite{bickel2015mathematical}). In parametric settings, one can invert hypothesis tests to obtain confidence intervals. Classical theory shows that for a fixed sample size and a prescribed Type I error, there exist UMP (uniformly most powerful) tests that maximize power, often restricting attention to particular families (e.g., such as unbiased tests in exponential families). This optimality, when translated to CIs, is usually framed as minimizing expected width, subject to some restriction on the class of CI, such as unbiasedness (See Definition 9.3.7 in \cite{casella2024statistical}). Further, Classic theory indicates that no CI can uniformly minimize random width at fixed $N$ and $\delta$ (see page 250 in \cite{bickel2015mathematical}). To overcome this negative result, we use asymptotic analysis. Our framework adopts a stability assumption, ensuring that the interval endpoints converge to deterministic values as $ \lim_{\delta \to 0} N_\delta \to \infty$. This makes the limiting width a well-defined metric and allows us to show that no stable policy can outperform our KL-based intervals, asymptotically. In the context of inverting a test statistic, existing results often allow for thresholds that depend intricately on sample size or distributional assumptions. The classical literature, for the most part, focuses on parametric families. In contrast, we show that a uniform threshold can be used in KL-divergence-based constructions for both parametric and nonparametric families with bounded support, yielding asymptotically optimal intervals. 

Additionally, there is a vast amount of literature on constructing CIs that proceeds by inverting finite-sample concentration inequalities. For the probability distributions with bounded support, one can use various concentration inequalities such as the Hoeffding and Bernstein inequalities (see \cite{hoeffding1994probability}, \cite{bennett1962probability}, \cite{boucheron2009bernstein}, \cite{audibert2009exploration} and \cite{maurer2009empirical}). In general, using the Chernoff bound, one can construct a CI of the mean if the upper bound on MGF is known. For the heavy tail distributions, one can get a CI of the mean using Markov's inequality and Chebyshev's inequality (see \cite{catoni2012challenging} and \cite{chen2021generalized}). Bootstrap methods \citep{efron1994introduction} are a popular approach 
for constructing confidence intervals, but their coverage guarantees are 
only asymptotically valid as $N \to \infty$. In our setting, we require exact 
finite-sample guarantees for any fixed $N$ and $\delta$, which standard bootstrap 
methods cannot provide.

Recently, \cite{waudby2024estimating} transform a CI construction problem into a betting problem for the wealth maximization process for bounded probability distributions. A recent work \cite{gupta2023minimax} utilizes a different notion of optimality (mini-max)  for the width of CI in location parameter families. We compare our policy's performance with some of these relevant existing methods in Section \ref{sec_num}.

There appears to be only one work in the literature that aims to characterize a lower bound on the width of a CI in terms of a distribution-dependent complexity term: \cite{shekhar2023near}. As noted in Remark 4.4 in \cite{shekhar2023near}: “To the best of our knowledge, this is the first result providing an explicit characterization of the smallest achievable width of CI in terms of a distribution-dependent complexity term.”

However, their bound is loose (in the sense that we strengthen it by a multiplicative constant asymptotically) and relies on the strong assumption that the width of the CI is deterministically bounded by a specific function for any given sample size $N$. In contrast, our lower bound on the limiting width is asymptotic but \textbf{remains tight} (in the sense the proposed CI construction policy matches the lower bound asymptotically), under a much milder assumption concerning the stability of the policies (see Section \ref{sec_num}).

The paper is organized as follows: In the next section, we formally introduce the setup and notation and focus on the single-parameter exponential family. In Section \ref{sec_lower_bound}, we present our results on the minimum limiting width for CIs in the three learning regimes: no learning, sufficient learning, and complete learning. In Section \ref{sec_kl}, we provide a recipe for constructing CIs using KL divergence-based concentration inequalities (see \eqref{CI_concen}) that achieve the minimum limiting width in the sufficient learning regime and zero limiting width in the complete learning regime. In Section \ref{sec_random}, we extend our results to the setting when there is a random cost associated with the sample collection. In Section \ref{sec_generalization}, we extend our results to the non-parametric distribution families with bounded support and known bound on $(1+\epsilon)^{\rm th}$ moment. We also discuss the asymptotic optimality of KL divergences-based CI construction policies in these settings. In Section \ref{sec_num}, we present our numerical studies. In Section \ref{sec_conclusion}, we present our conclusions and outline future research directions. In Appendix \ref{one_sided_CI}, we extend our results to the setting when we want to construct a one-sided CI of the mean.

\section{SETUP AND NOTATION}\label{sect_set_up}
Let \( X_1, X_2, \dots \) be i.i.d. copies of a random variable \( X \) with distribution \( \nu \) and mean \( m(\nu) = \mathbb{E}_{\nu}[X] = \mu \). For $ n \geq 1$, let ${\cal F}_n$ denote the information contained in the $\sigma$-algebra generated by  $\{X_k,  k \leq n\}$. We now state a formal definition of confidence intervals.
\subsection{Description of the policy space for construction of CI}
Let $\Pi_{\rm CI}$ denote the collection of methods/policies for constructing CIs. Let $[\widehat{\mu}^{\pi}_L(N,\delta), \widehat{\mu}^{\pi}_R(N,\delta)]$ denote the estimated CI after observing $X_1,X_2,\ldots,X_N$ under a policy $\pi$ for any given $\delta \in (0,1)$. For any policy $\pi \in \Pi_{\rm CI}$ must satisfy the following for any given $\delta \in (0,1)$,
$$ \forall n \in \mathbb{N} : \mathbb{P}_{\nu}(\mu \in [\widehat{\mu}^{\pi}_L(n,\delta), \widehat{\mu}^{\pi}_R(n,\delta)]) \geq 1-\delta.$$

Here, $\mathbb{P}_{\nu}(\cdot)$ denotes the probability measure induced by the environment $\nu$. Importantly, the above coverage requirement is non-asymptotic: it must hold for every fixed $N$ and $\delta$.

Recall that our objective is to characterize the width of CI, i.e., $ \widehat{\mu}^{\pi}_R(N,\delta) -  \widehat{\mu}^{\pi}_L(N,\delta)$. As stated in the Literature review section, for a given $N$ and $\delta$, a tight characterization of the width of CI is analytically intractable, and thus we consider the asymptotic regime where $N \to \infty$ as $\delta \to 0$ and aim to characterize the limiting width of CIs. Henceforth, we denote $N$ as $N_\delta$. To start the analysis, we first make a stability assumption over the space of policies which construct CIs, enabling us to derive a lower bound on the limiting width. Under this assumption, the limiting width of CIs becomes a deterministic constant for any CI construction policy. This assumption requires that for a given environment $\nu$, the boundaries of the CI, $\widehat{\mu}^{\pi}_L(N_\delta,\delta), \widehat{\mu}^{\pi}_R(N_{\delta},\delta)$ converge in probability to deterministic points as $\delta$ approaches $0$ and $N_\delta$ approaches infinity. In the Appendix, we show that many popular policies for constructing CI are stable.

\begin{definition}[Stability] \label{def_stability}
Let $N_\delta \to \infty$ as $\delta \to 0$. For a given distribution $\nu$ with
mean $\mu$, a policy $\pi \in \Pi_{\rm CI}$ is called \textbf{stable} if the CI it
constructs, denoted by $[\widehat{\mu}^{\pi}_L(N_\delta, \delta),\,
\widehat{\mu}^{\pi}_R(N_\delta, \delta)]$, satisfies: if
\[
    \lim_{\delta \to 0} \frac{N_\delta}{\log(1/\delta)} \in [0, \infty],
\]
then
\[
    \widehat{\mu}^{\pi}_L(N_\delta, \delta) \overset{p}{\longrightarrow} \mu_L^{\pi}(\nu)
    \qquad \text{and} \qquad
    \widehat{\mu}^{\pi}_R(N_\delta, \delta) \overset{p}{\longrightarrow} \mu_R^{\pi}(\nu),
\]
where $\mu_L^{\pi}(\mu) \leq \mu$ and $\mu_R^{\pi}(\mu) \geq \mu$ are deterministic
constants (with $-\infty$ and $\infty$ included).
\end{definition}

We denote the collection of policies in the set $\Pi_{\rm CI}$ which are stable as $\Pi^{s}_{\rm CI}$. It is worth noting that $\mu_{R}^{\pi}(\nu) -\mu_{L}^{\pi}(\nu)$ denotes the limiting width of CI for $\pi \in \Pi^{s}_{\rm CI} $ in the asymptotic regime where the sample size, i.e, $N_\delta$, scales to $\infty$ as $\delta \to 0$. We now assume that $\nu$ belongs to the canonical single-parameter exponential family $\mathbf{S}$. In Section \ref{sec_generalization}, we later generalize our results to non-parametric distributions.

Specifically, $\mathbf{S}$ is defined as: $   \mathbf{S} = \left\{ p_{\theta} : \theta \in \Theta \subseteq \mathbb{R}, \frac{dp_{\theta}(x)}{d\xi} = \exp(\theta \cdot x - b(\theta)) \right\}$, where $\xi$ is a fixed reference measure on $\mathbb{R}$, and $b(\theta)$ is a known, twice-differentiable, strictly convex function. The set $\Theta$ is:
$
\left\{ \theta \in \mathbb{R}:  \int_{\mathbb{R}} |x| \exp(\theta \cdot x) d\xi(x) < \infty \right\}.$ This condition ensures that $p_{\theta}$ forms a well-defined probability distribution with a finite mean. We assume $\Theta = (\underline{\theta}, \overline{\theta})$ is an open interval. Common distributions in this family include Bernoulli, Poisson, Gaussian (with known variance), and Gamma (with a known shape parameter) (see \cite{cappe2013kullback} for more details on single-parameter exponential family). For $\nu = p_{\theta} \in \mathbf{S}$, the mean $\mu = \mathbb{E}_{\nu}[X]$ is unknown. It is known that $\mu(\theta) = b'(\theta)$, which, due to the strict convexity of $b(\theta)$, is a strictly increasing function. This allows us to define the inverse function $\theta(\mu)$. We assume that the support of the mean is denoted by $\mathbf{O} = (\underline{\mu}, \overline{\mu})$.

\noindent\textbf{Divergence function:} Let $KL(p_{\theta}, p_{\tilde{\theta}})$ denote the KL divergence. We define:

\begin{align*}
d(\mu, \tilde{\mu}) 
&= KL(p_{\theta(\mu)}, p_{\theta(\tilde{\mu})}) \\
&= b(\theta(\tilde{\mu})) - b(\theta(\mu))
   - b'(\theta(\mu)) \, (\theta(\tilde{\mu}) - \theta(\mu)).
\end{align*}

for $\mu, \tilde{\mu} \in \mathbf{O}$. Key properties of $d(\mu, \tilde{\mu})$ include its strict quasi-convexity in the second argument and the fact that $d(\mu, \tilde{\mu}) > 0$ for all $\tilde{\mu} \neq \mu$, with $d(\mu, \mu) = 0$. It is worth noting that, since for $\nu \in \mathbf{S}$, the distribution is uniquely determined by its mean $\mu$. Therefore, for simplicity, we denote the limiting CI of a $\pi \in \Pi^{s}_{\rm CI}$, $[\mu_L^{\pi}(\nu), \mu_R^{\pi}(\nu)]$, as $[\mu_L^{\pi}(\mu), \mu_R^{\pi}(\mu)]$.



 
\section{MAIN RESULTS ON THE MINIMUM LIMITING WIDTH OF CONFIDENCE INTERVAL IN DIFFERENT REGIMES}\label{sec_lower_bound}

Recall that \( \nu \), with mean \( \mu \), represents the true distribution from which samples are generated. In other words, \( \nu \) is the true underlying but unknown environment. To start the analysis of the minimum limiting width, we define an alternate environment $\tilde{\nu}$ such that $\mathbb{E}_{\tilde{\nu}}[X] = \tilde{\mu} \neq \mu$. Using the data processing inequality (see \cite{cover1991elements}), it follows that for a given $\delta \in (0,1)$ and any alternate environment $\tilde{\nu}$ and any event $\mathcal{E}_\delta$, we have 
	\begin{equation} \label{eqn_lower:01}
N_\delta \cdot d(\mu, \tilde{\mu})  \geq \sup_{\mathcal{E}_\delta \in \mathcal{F}_{N_\delta}} \phi(\mathbb{P}_{\nu}(\mathcal{E}_\delta),
	\mathbb{P}_{\tilde{\nu}}(\mathcal{E}_\delta)),
  	\end{equation}
where  $\phi(p_1,p_2) \triangleq  p_1 \log\frac{p_1}{p_2} +(1-p_1)\log \left(\frac{1-p
_1}{1-p_2}\right)$ for $p_1, p_2 \in (0,1)$. To utilize the above result, consider any policy $\pi \in \Pi^{s}_{\rm CI}$. We now define the set of alternate environments $K(\mu_L^\pi(\mu),\mu_R^\pi(\mu)) = K_1(\mu_L^\pi(\mu)) \cup K_2(\mu_R^\pi(\mu))$, where $  K_1(\mu_L^\pi(\mu)) = \{ \tilde{\nu}: \tilde{\nu} \in \mathbf{S}, \tilde{\mu} < \mu_L^\pi(\mu)\}$ and  $K_2(\mu_R^\pi(\mu)) = \{ \tilde{\nu}: \tilde{\nu} \in \mathbf{S}, \tilde{\mu} > \mu_R^\pi(\mu)\}$. We utilize \eqref{eqn_lower:01} for $\tilde{\nu} \in K(\mu_L^\pi(\mu),\mu_R^\pi(\mu))$ and $\mathcal{E}_\delta = \{\tilde{\mu} \notin [\widehat{\mu}^{\pi}_L(N_\delta,\delta), \widehat{\mu}^{\pi}_R(N_\delta,\delta)]\}$, where recall that $\tilde{\mu}= m(\tilde{\nu})$. Using the fact that $\pi \in \Pi^{s}_{\rm CI}$, we get that  $\mathbb{P}_{\tilde{\nu}}(\mathcal{E}_\delta)\leq \delta$. Since $\tilde{\nu} \in K(\mu_L^\pi(\mu),\mu_R^\pi(\mu))$, it implies that $\tilde{\mu} >  \mu_R^\pi(\mu)$ or $ \tilde{\mu}  < \mu_L^\pi(\mu)$. Further, as $\pi$ is a stable policy and hence it implies $\mathbb{P}_{\nu}(\mathcal{E}_\delta) \approx 1 $ for small $\delta$. Using \eqref{eqn_lower:01} and dividing both sides with $\log(1/\delta)$ and taking $\delta \to 0$, we get different lower bounds on the limiting width depending upon the scaling of $N_\delta$. We now present the formal result and the rigorous proof is given in the Appendix \ref{appen_lower_bound}.

\begin{theorem} \label{thm_spef_lower_bound_width}
    For a given $\nu \in \mathbf{S}$ with mean $\mu$, and any $\pi \in \Pi^{s}_{\rm CI}$, the following holds:

    a) \textbf{ No learning regime :} If $\lim_{\delta \to 0} \frac{N_\delta} {\log(1/\delta)}\to 0$ then, $   \left[ \mu_{R}^{\pi}(\mu) - \mu_{L}^{\pi}(\mu) \right] =  \overline{\mu} - \underline{\mu}.$ Further, $\lim_{\delta \to 0 } \widehat{\mu}^{\pi}_L(N_\delta,\delta) \overset{p}{\to} \underline{\mu} $ and $ \lim_{\delta \to 0 } \widehat{\mu}^{\pi}_R(N_\delta,\delta) \overset{p}{\to} \overline{\mu}.$

   b) \textbf{ Sufficient learning regime :} If $\lim_{\delta \to 0} \frac{N_\delta}{\log(1/\delta)} \to k$ for $k \in (0,\infty)$, then we have, 
  \begin{equation} \label{main_sup_1}
    \left[ \mu_{R}^{\pi}(\mu) - \mu_{L}^{\pi}(\mu) \right] \geq   \mu_{R}^{*}(\mu,k) - \mu_L^*(\mu, k),
  \end{equation}
      where, $ \mu_{L}^{*}(\mu,k) <\mu $ and $ \mu_{R}^{*}(\mu,k) > \mu$ uniquely solve the following system of equations,
\begin{equation} \label{eqn_kl_symmetric}
    d(\mu,\mu_{R}^{*}(\mu,k) ) = d(\mu, \mu_{L}^{*}(\mu,k) )   =\frac{1}{k}.
\end{equation}
\end{theorem}
For the case when the sample size scales at the rate of $\log(1/\delta)$, as $\delta \to 0$, we obtain a lower bound on the limiting width of the CI, given by $\mu_{R}^{*}(\mu,k) - \mu_{L}^{*}(\mu,k)$. It follows from the quasi-convexity of $d(\mu,x)$ in $x$, implying that $\mu_{R}^{*}(\mu,k) - \mu_{L}^{*}(\mu,k)$ decreases as $k$ increases. In the proof, we first demonstrate that $\mu_{R}^{\pi}(\mu) \geq \mu_{R}^{*}(\mu,k)$ for any $\pi \in \Pi^{s}_{\rm CI}$. We then show that $\mu_{L}^{\pi}(\mu) \leq \mu_{L}^{*}(\mu,k)$ for any $\pi \in \Pi^{s}_{\rm CI}$. The key idea of the proof is that for $\mu_{L}^{\pi}(\mu) > \mu_{L}^{*}(\mu,k)$ and $\mu_{R}^{\pi}(\mu) < \mu_{R}^{*}(\mu,k)$, we get a contradiction with \eqref{eqn_lower:01}. Hence, $[\mu_{L}^{*}(\mu,k),\mu_{R}^{*}(\mu,k)]$ can be interpreted as a subset of the CI constructed by any policy $\pi \in \Pi^{s}_{\rm CI}$ in the sufficient learning regime. Later, we show that our proposed policy $\pi_1$ has $[\mu_{L}^{*}(\mu,k), \mu_{R}^{*}(\mu,k)]$ as the limiting CI in the sufficient learning regime, proving that the above lower bound on the limiting width is tight.

Our main novelty and contribution in the lower bound proof lies in the introduction of the stability notion for CI construction methods, and in further leveraging this concept together with \eqref{eqn_lower:01} (data processing inequality). This stability concept arises very naturally in the asymptotic regime we study, which itself has not been explored in the context of confidence interval width. Importantly, the stability assumption is very mild, as it is trivially satisfied by standard concentration bound-based methods, such as those based on Hoeffding or empirical Bernstein inequalities (see Appendix).

Further, due to the stability notion, we obtain a tighter lower bound than that of Proposition 4.3 of \cite{shekhar2023near}. To be precise, the authors in [4] derive one-sided lower bounds on the width using a similar argument to the data processing inequality, and then take the maximum of the left and right deviations. Ideally, one should be able to combine the lower bounds on the width from the left and right deviations, but unfortunately, it is not possible trivially. However, the stability assumption allows us to consider joint elimination of hypotheses on both sides, leading to a tighter bound that equals the total limiting width. This results in a factor of two improvement in the Gaussian case (as shown in Section \ref{sec_num}). Thus, while our proof begins with a known inequality, the strength of our result comes from this careful asymptotic characterization using stability, which leads to tighter bounds.

\subsection{Complete learning regime}

 It is worth noting that for the case when $\lim_{\delta \to 0} \frac{N_\delta} {\log(1/\delta)} \to \infty$, similar to the proof of above theorem, we get a trivial lower bound on the limiting width, i.e.,  $    \left[ \mu_{R}^{\pi}(\mu) - \mu_{L}^{\pi}(\mu) \right] \geq 0$ for any $\pi\in \Pi^{s}_{\rm CI}$. Further, our proposed policy $\pi_1$ (see Section \ref{sec_kl}) has zero limiting width in this regime. Hence, we denote this regime as the complete learning regime. Additionally, in the complete learning regime, we characterize the rate at which the CI width converges to zero. Under certain technical assumptions on CI construction policies, we establish the fastest achievable convergence rate. Furthermore, we demonstrate that the CI width under our proposed policy $\pi_1$ attains this optimal rate as it approaches zero (see Appendix \ref{sub_sec_rate_analysis} for formal results).

\begin{remark}
It is worth noting that, for any theoretical optimality guarantee for a CI procedure, one can think of four
distinct formulations: (i)~fixed $N$ and $\delta$; (ii)~fixed $N$ and $\delta \to 0$;
(iii)~fixed $\delta$ and $N \to \infty$; and (iv)~$N \to \infty$ and $\delta \to 0$
jointly. Since classical theory shows that no CI can uniformly minimize random width for
    fixed $N$ and $\delta$ (see Section \ref{sec_lit}), formulation~(i) does not yield
meaningful results. Formulations~(ii) and~(iii) are relatively trivial: if only
$\delta \to 0$ while $N$ remains fixed, the results fall into the no-learning regime;
if only $N \to \infty$ while $\delta$ remains fixed, the results correspond to the
complete learning regime, where the width shrinks to zero. The substantive and non-trivial setting is therefore formulation~(iv), where $N \to \infty$ and
$\delta \to 0$ jointly. We find that if
$   \lim_{\delta \to 0} \frac{N_\delta}{\log(1/\delta)} = k \in [0, \infty],$
the three different regimes, no learning ($k = 0$), sufficient learning
($k \in (0,\infty)$), and complete learning ($k = \infty$), arise naturally, and our
results precisely characterize the minimum achievable limiting CI width in each case.
\end{remark}

At last, we discuss the possible connections of our lower bound result with the classical Cram\'er--Rao lower bound in Appendix \ref{appen_cramer}

\section{ ASYMPTOTIC OPTIMALITY OF KL DIVERGENCE BASED CONSTRUCTION OF CI: DESCRIPTION OF METHOD $\pi_1$}  \label{sec_kl}
In this section, we describe the method/policy $\pi_1$. It is well known that for $\nu \in \mathbf{S}$, the sample average is the MLE of the mean, $\hat{\mu}_n = \frac{\sum_{t=1}^{n} X_t}{n}.$ For $\nu \in \mathbf{S}$, we utilize the concentration inequality based on KL divergences:
\begin{equation} \label{CI_concen}
    \mathbb{P}_{\nu} ( n \,d(\hat{\mu}_n,\mu) \geq \beta(\delta)) \leq \delta.
\end{equation}
Here $\beta(\delta) = \log(2/\delta)$ is a well-chosen function so the above holds. This can be derived from Lemma 4 in \cite{menard2017minimax} or Theorem 4 in \cite{busa2019optimal} and the proof is based on applying Markov's inequality and the structure of $d(\cdot,\cdot)$ function.
To see, how we construct CI from the above concentration inequality, we formally define $\mu^{\pi_1}_{L} (n,\delta)$ and  $\mu^{\pi_1}_{R} (n,\delta)$ as follows: 
\begin{equation} \label{policy_pi_1}
    \mu^{\pi_1}_{R} (n,\delta) \triangleq \max \{q > \hat{\mu}_n :n d(\hat{\mu}_n , q) \leq \beta(\delta) \} \textrm{ and }
\end{equation}
$$ \mu^{\pi_1}_{L} (n,\delta) \triangleq \min \{q < \hat{\mu}_n :nd(\hat{\mu}_n , q) \leq \beta( \delta) \}.$$


Now we state the formal result related to the limiting width of the CI construction under policy $\pi_1$. 

\begin{theorem} \label{thm_kl}
  The policy $\pi_1$ has following properties:
    
    a) $\pi_1 \in \Pi^{s}_{\rm CI}$.
    
    b) If $\lim_{\delta \to 0} \frac{N_\delta} {\log(1/\delta)}\to k$ for $k \in (0,\infty)$, then we have, $   \mu_{R}^{\pi_1}(\mu)  =  \mu_{R}^{*}(\mu,k) $ and $ \mu_{L}^{\pi_1}(\mu) =  \mu_{L}^{*}(\mu,k)$, where, $ \mu_{L}^{*}(\mu,k) <\mu $ and $ \mu_{R}^{*}(\mu,k) > \mu$ uniquely solve \eqref{eqn_kl_symmetric}.

c) If $\lim_{\delta \to 0} \frac {N_\delta} {\log(1/\delta)}\to \infty$, then we have, 
$  \mu_{R}^{\pi_1}(\mu) - \mu_{L}^{\pi_1}(\mu)  = 0$.
\end{theorem}

The above theorem shows that \( \pi_1 \) achieves the minimum limiting width in the sufficient learning regime and zero limiting width in the complete learning regime, making it an asymptotically optimal policy for CI construction of the mean. Next, we study the setting where samples are costly.

\section{RANDOM SAMPLING COST} \label{sec_random}

Motivated by the applications discussed in Section \ref{sec:intro}, this section focuses on the problem of constructing a CI for the mean when there is a random cost associated with obtaining samples. In this setting, instead of a fixed sample size \( N \), we assume a cost budget of \( C \) units. As in the previous setting, we will scale \( C \) with \( \delta \), and henceforth denote the cost budget as \( C_\delta \). Each sample incurs a random cost $c_i$, where $c_i$ for $i = 1, 2, 3, \dots$ are i.i.d. copies from an unknown cost distribution $\mathcal{C}$ with positive support and mean $\overline{c} > 0$. We assume that the distributions \( \nu \) and \( \mathcal{C} \) are independent. Consequently, \( X_i \) and \( c_i \) for \( i = 1, 2, 3, \dots \) are also independent.

Our goal is to establish the lower bound on the limiting width of any CI construction policy that ensures the mean lies within the interval with probability \( (1 - \delta) \) given the budget \( C_\delta \). To do this, we define a random time, interpreted as the maximum number of samples that can be collected within the given budget \( C_\delta \): $\tau_\delta = \sup \{ n \in \mathbb{Z}^{+}: \sum_{i = 1}^{n}c_i \leq C_\delta \}.$

Let $[\widehat{\mu}^{\pi}_L(\tau_\delta,\delta), \widehat{\mu}^{\pi}_R(\tau_\delta,\delta)]$ denote the estimated CI after observing $X_1,X_2,\ldots,X_{\tau_\delta}$ under a policy $\pi$ for any given $\delta \in (0,1)$. Let $\hat{\Pi}_{\rm CI}$ denote the collection of CI construction policies such that the following holds for any $\pi \in \hat{\Pi}_{\rm CI}$ and for any given $\delta \in (0,1)$, $\mathbb{P}_{\nu}(\mu \in [\widehat{\mu}^{\pi}_L(\tau_\delta,\delta), \widehat{\mu}^{\pi}_R(\tau_\delta,\delta)]) \geq 1-\delta.$
    \begin{definition} \label{def_stability}
   Let $C_\delta \to \infty$ as $\delta \to 0$. For a given distribution $\nu$ with mean $\mu$, for any $\pi \in \hat{\Pi}_{\rm CI}$, $\pi$ is called a stable policy if following holds,
   
 $\lim_{\delta \to 0 } \widehat{\mu}^{\pi}_L(\tau_\delta,\delta) \overset{p}{\to} \mu_{L}^{\pi}(\nu),$
    $\lim_{\delta \to 0 } \widehat{\mu}^{\pi}_R(\tau_\delta,\delta) \overset{p}{\to} \mu_{R}^{\pi}(\nu),$
where $\mu_{L}^{\pi}(\nu) \leq  \mu$ and   $\mu_{R}^{\pi}(\nu) \geq  \mu$ are constants ($\infty$ and $-\infty$ included).
\end{definition}

We denote the collection of policies in the set $\hat{\Pi}_{\rm CI}$ which are stable as $\hat{\Pi}^{s}_{\rm CI}$. It is worth noting that $\mu_{R}^{\pi}(\nu) -\mu_{L}^{\pi}(\nu)$ denotes the limiting width of CI for $\pi \in \Pi^{s}_{\rm CI} $ in the asymptotic regime where the sample size, i.e, $C_\delta \to \infty$ as $\delta \to 0$. We again assume that $\nu$ belongs to the canonical single-parameter exponential family $\mathbf{S}$. Hence, for simplicity, we denote 
$ \mu_{L}^{\pi}(\nu)$ and $ \mu_{R}^{\pi}(\nu)$ as $ \mu_{L}^{\pi}(\mu)$ and $ \mu_{R}^{\pi}(\mu)$, respectively. In Section \ref{sec_generalization}, we generalize our results for non-parametric distributions.
Now we state our key result on the lower bound on the limiting width of any $\pi \in \hat{\Pi}^{s}_{\rm CI}$.
\begin{theorem} \label{thm_spef_lower_bound_width_cost}
    For a given $\nu \in \mathbf{S}$ with mean $\mu$, a cost distribution $\mathcal{C}$ with mean $ 0<\overline{c} < \infty$, and any $\pi \in \hat{\Pi}^{s}_{\rm CI}$, the following holds:

    a) \textbf{ No learning regime :} If $\lim_{\delta \to 0} \frac{C_\delta} {\log(1/\delta)}\to 0$, then $  \left[ \mu_{R}^{\pi}(\mu) - \mu_{L}^{\pi}(\mu) \right] =  \overline{\mu} - \underline{\mu}.$ Further, $\lim_{\delta \to 0 } \widehat{\mu}^{\pi}_L(\tau_\delta,\delta) \overset{p}{\to} \underline{\mu} $ and $ \lim_{\delta \to 0 } \widehat{\mu}^{\pi}_R(\tau_\delta,\delta) \overset{p}{\to} \overline{\mu}.$

   b) \textbf{ Sufficient learning regime :} If $\lim_{\delta \to 0} \frac{C_\delta}{\log(1/\delta)} \to k$ for $k \in (0,\infty)$, then
  \begin{equation} \label{main_sup_1_cost}
    \left[ \mu_{R}^{\pi}(\mu) - \mu_{L}^{\pi}(\mu) \right] \geq   \mu_{R}^{*}(\mu,k, \overline{c}) - \mu_L^*(\mu, k, \overline{c}),
  \end{equation}
      where, $ \mu_{L}^{*}(\mu,k, \overline{c}) <\mu $ and $ \mu_{R}^{*}(\mu,k, \overline{c}) > \mu$ uniquely solve the following system of equations,
\begin{equation} \label{eqn_kl_symmetric_cost}
    d(\mu,\mu_{R}^{*}(\mu,k, \overline{c}) ) = d(\mu, \mu_{L}^{*}(\mu,k, \overline{c}) )   =\frac{\overline{c}}{k}.
\end{equation}
\end{theorem}

The presence of an additional average sample collection cost, $\overline{c}$, in \eqref{eqn_kl_symmetric_cost} differentiates the above result from Theorem \ref{thm_spef_lower_bound_width}. Furthermore, the above result indicates that the limiting width of the confidence interval is invariant to the distribution of the cost and only depends on its mean. Apart from the ideas in the proof of Theorem \ref{thm_spef_lower_bound_width}, three additional key ideas are required to prove the above theorem. The first is the use of renewal theory to study the scaling of $\tau_\delta$ with $C_\delta$. The second is the fact that \eqref{eqn_lower:01} holds for a stopping time, and in this setting, $\tau_\delta + 1$ is also a stopping time. The last idea is to extend \eqref{eqn_lower:01} to the joint distribution of $\nu$ and $\mathcal{C}$.

\begin{remark}
Analogous to Remark 1, here when $\lim_{\delta \to 0} \frac{C_\delta} {\log(1/\delta)} \to \infty$, we have a trivial lower bound on the limiting width, i.e.,  $    \left[ \mu_{R}^{\pi}(\mu) - \mu_{L}^{\pi}(\mu) \right] \geq 0$ for any $\pi\in \hat{\Pi}^{s}_{\rm CI}$. Further, our modified policy $\hat{\pi}_1$ (see below) has zero limiting width in this regime. Again, we denote this regime as the complete learning regime. 
\end{remark}

 In this extended setting, we propose a modified method, denoted as $\hat{\pi}_1$, which is identical to $\pi_1$, except with a modified $\beta(n, \delta)$ replacing $\beta(\delta)$ in \eqref{policy_pi_1}. This modification is necessary for the following reason: the number of samples is random in this setting, and the concentration inequality used for a fixed number of samples, as given in \eqref{CI_concen}, is no longer valid. To address this, we employ an anytime-valid concentration inequality. Specifically, we use $\beta(n, \delta)$ from (14) in \cite{kaufmann2021mixture}, defined as:
$\beta(n, \delta) = 3\log(1+\log(n)) + \mathcal{T}(\log(1/\delta)).
$

Here, the function $\mathcal{T}(x): \mathbb{R}^+ \to \mathbb{R}^+$ is given by:$
\mathcal{T}(x) = 2\Tilde{\psi}_{3/2}\left(\frac{x + \log(2\zeta(2))}{2}\right),
$
$\text{where } \zeta(2) = \sum_{n=1}^{\infty} n^{-2}$ and, 

\[
\Tilde{\psi}_y(x) = 
\begin{cases}
e^{1/\psi^{-1}(x)}\psi^{-1}(x) & \text{if } x \geq \psi^{-1}(1/\ln{y}), \\
y(x - \ln{\ln{y}}) & \text{otherwise},
\end{cases}
\]

for any \( y \in [1, e] \) and \( x \geq 0 \). The function \( \psi(u) = u - \ln{u} \) has an inverse \( u = \psi^{-1}(z) \) for \( z \geq 1 \). As shown by \cite{kaufmann2021mixture}, it holds that \( \lim_{x \to \infty} \frac{\mathcal{T}(x)}{x} = 1 \). Now we state the formal result that the limiting width of the CI under policy $\hat{\pi}_1$ matches the lower bound.

\begin{theorem} \label{thm_kl_cost}
  The policy $\hat{\pi}_1$ has following properties:
    
    a) $\hat{\pi}_1 \in \hat{\Pi}^{s}_{\rm CI}$.
    
    b) If $\lim_{\delta \to 0} \frac{C_\delta} {\log(1/\delta)}\to k$ for $k \in (0,\infty)$, then we have, $   \mu_{R}^{\hat{\pi}_1}(\mu) =  \mu_R^*(\mu, k,\overline{c}),$ and $  \mu_{L}^{\hat{\pi}_1}(\mu)  =  \mu_{L}^{*}(\mu,k,\overline{c}) $ where, $ \mu_{L}^{*}(\mu,k, \overline{c}) <\mu $ and $ \mu_{R}^{*}(\mu,k,\overline{c}) > \mu$ uniquely solve \eqref{eqn_kl_symmetric_cost}.

c) If $\lim_{\delta \to 0} \frac {C_\delta} {\log(1/\delta)}\to 0$, then we have, 
$  \mu_{R}^{\hat{\pi}_1}(\mu) - \mu_{L}^{\hat{\pi}_1}(\mu)  = 0$.
\end{theorem}

\section{GENERALIZATION TO NON-PARAMETRIC DISTRIBUTIONS} \label{sec_generalization}

In this section, we generalize our results for the case when the underlying distribution belongs to a non-parametric family. In particular, we consider two such settings. The first one is the family of probability distributions with bounded support in $[0,1]$, we denote it as $\mathbf{B}$. The second one is the family of probability distributions with $(1+\epsilon)^{\rm th}$ moment bounded and the bound is known and given by a constant $\Gamma$. We denote this family as $\mathbf{H}$. It is worth noting that $KL_{\rm inf}(\nu,\mathbf{B},x)$ is a well-studied in \cite{honda2010asymptotically} and \cite{jourdan2022top} for $\nu \in \mathbf{B}$. Further, \cite{agrawal2022bandits} studies $KL_{\rm inf}(\nu,\mathbf{H},x)$ for $\nu \in \mathbf{H}$. A primer on $KL_{\mathrm{inf}}(\nu,\mathbf{B},x)$ and $KL_{\mathrm{inf}}(\nu,\mathbf{H},x)$ and how they are computed is provided in the Appendix \ref{appen_bounded},\ref{appen_heavy}. The analysis in this section is similar to the setting where $\nu \in \mathbf{S}$, with $d(\mu, x)$ now replaced by $KL_{\rm inf}(\nu, \mathbf{P}, x)$ (defined below). 

Given a probability distribution family $\mathbf{P}$, an outcome distribution $\nu \in \mathbf{P}$ and $x \in \mathbf{R}$, let, 
\begin{equation} \label{eq_d_inf}
KL_{\rm inf} (\nu,\mathbf{P},x) =
\begin{cases} 
    \inf\limits_{\substack{\kappa \in \mathbf{P} \\ m(\kappa)\geq x}} KL (\nu,\kappa), & \text{if } x \geq m(\nu) \\[4pt]
    \inf\limits_{\substack{\kappa \in \mathbf{P} \\ m(\kappa)\leq x}} KL (\nu,\kappa), & \text{if } x < m(\nu).
\end{cases}
\end{equation}

$KL_{\rm inf} (\nu,\mathbf{P},x)$ is the minimum amongst the KL divergences between a given distribution $\nu$ and all distributions in the same family $\mathbf{P}$ which have higher mean than $x$ if $x \geq m(\nu)$. It is similarly defined for $x < m(\nu)$.
\subsection{Results on lower bound on limiting width of CI}\label{sub_sec_general_lower_bound}

For the setting, when $\nu \in \{\mathbf{B}, \mathbf{H}\}$, Theorem \ref{thm_spef_lower_bound_width_cost} extends as follows.

\begin{theorem} \label{thm_generalization_lower}
    For a given $\nu \in \{\mathbf{B}, \mathbf{H}\}$ with mean $m(\nu) = \mu$, a cost distribution $\mathcal{C}$ with mean $ 0<\overline{c} < \infty$, and any $\pi \in \hat{\Pi}^{s}_{\rm CI}$, the following holds:

    a) \textbf{ No learning regime :} If $\lim_{\delta \to 0} \frac{C_\delta} {\log(1/\delta)}\to 0$, then $\left[ \mu_{R}^{\pi}(\nu) - \mu_{L}^{\pi}(\nu) \right] =  \overline{\mu} - \underline{\mu}.$

    Further, $\lim_{\delta \to 0 } \widehat{\mu}^{\pi}_L(\tau_\delta,\delta) \overset{p}{\to} \underline{\mu} $ and $ \lim_{\delta \to 0 } \widehat{\mu}^{\pi}_R(\tau_\delta,\delta) \overset{p}{\to} \overline{\mu}.$

   b) \textbf{ Sufficient learning regime :} If $\lim_{\delta \to 0} \frac{C_\delta}{\log(1/\delta)} \to k$ for $k \in (0,\infty)$, then
  \begin{equation} \label{main_sup_1_cost}
    \left[ \mu_{R}^{\pi}(\nu) - \mu_{L}^{\pi}(\nu) \right] \geq   \mu_{R}^{*}(\nu,k, \overline{c}) - \mu_L^*(\nu, k, \overline{c}),
  \end{equation} $ \mu_{L}^{*}(\nu,k,\overline{c}) <\mu $ and $ \mu_{R}^{*}(\nu,k,\overline{c}) > \mu$ uniquely solve the following system of equations,
\begin{equation} \label{eqn_kl_symmetric_2}
    KL_{\rm inf}(\nu,\mathbf{P},\mu_{R}^{*}(\nu,k,\overline{c}) ) = KL_{\rm inf}(\nu,\mathbf{P}, \mu_{L}^{*}(\nu,k,\overline{c}) )   =\frac{\overline{c}}{k},
\end{equation}
where, $\mathbf{P} = \mathbf{B}$ if $\nu \in \mathbf{B}$ and $\mathbf{P} = \mathbf{H}$ if $\nu \in \mathbf{H}$.
\end{theorem}

In this setting, when $\nu \in \{\mathbf{B}, \mathbf{H}\}$, we denote the limiting CI of a policy $\pi \in \hat{\Pi}_{\rm CI}^{s}$ as $[\mu_{L}^{\pi}(\nu), \mu_{R}^{\pi}(\nu)]$ as $\delta \to 0$ and $C_\delta \to \infty$.
\subsection{Asymptotic optimality of KLinf-based construction of CI} \label{sub_sec_general_kl_cs}

We use the following concentration inequality for the construction of the CI (see \cite{agrawal2022bandits}, \cite{orabona2023tight}, and \cite{jourdan2022top}), 
\begin{equation}\label{CS_conce_2}
    \mathbb{P}_{\nu} ( \exists \, n \in \mathbb{N}: \,n \, KL_{\rm inf}(\hat{\nu}_n, \mathbf{B},\mu) \geq \beta(n,\delta)) \leq \delta,
\end{equation}
 $\hat{\nu}_n$ denotes the empirical distribution after $n$ samples and $\beta(n,\delta) = 1+  \log\left(\frac{2(1+n)}{\delta}\right)$. Similar to the $\pi_1$ and $\hat{\pi}_1$, we now utilize the above KLinf-based concentration inequality to get a CI construction method denoted as $\pi_1^{\rm b}$. Let, $\mu^{\pi_1^{\rm b}}_{R} (n,\delta) \triangleq \max \{q > m(\hat{\nu}_n) : 
\quad n \, KL_{\rm inf}(\hat{\nu}_n, \mathbf{B}, q) \leq \beta(n,\delta) \}, $ and $\mu^{\pi_1^{\rm b}}_{L} (n,\delta) \triangleq \min \{q < m(\hat{\nu}_n) : 
\quad n \, KL_{\rm inf}(\hat{\nu}_n, \mathbf{B}, q) \leq \beta(n,\delta) \}.$ 
It follows that the reported CI is $[\mu^{\pi_1^{\rm b}}_{L} (n,\delta),\mu^{\pi_1^{\rm b}}_{R} (n,\delta)]$. An equivalent version of Theorem \ref{thm_kl} and Theorem \ref{thm_kl_cost} hold for our policy when $\nu \in \mathbf{B}$. The statements and proofs are given in the Appendix \ref{appen_generalization}. A similar CI construction method for $\nu \in \mathbf{H}$ holds using concentration bound in \cite{agrawal2022bandits}. Further, an equivalent version of Theorem \ref{thm_kl} and Theorem \ref{thm_kl_cost} hold for the policy $\pi_1^{\rm b}$ when $\nu \in \mathbf{H}$(with a minor technical assumption). See Appendix \ref{appen_generalization_1} for the description of the CI construction method, Theorem statements and proofs. A numerical study for the heavy tailed case, i.e., $\nu \in \mathbf{H}$ is provided in Appendix \ref{appen_num} as well.

\begin{remark}
    We provide a non-asymptotic analysis of the width of the CI for our policies under the cases where $\nu \in \mathbf{S}$ and $\mathbf{B}$ in Appendix \ref{appen_finite_confidence}.
\end{remark}
\section{NUMERICAL EXPERIMENTS} \label{sec_num}
Our numerical study has two objectives. First, to demonstrate numerically that our asymptotic lower bound is sharper and tighter than that of \cite{shekhar2023near} in our asymptotic regime. Second, to demonstrate the performance of the CI construction method $\pi_1$ and compare it with existing methods. 

\begin{figure}[ht]
    \centering
    \includegraphics[width=0.9\linewidth]{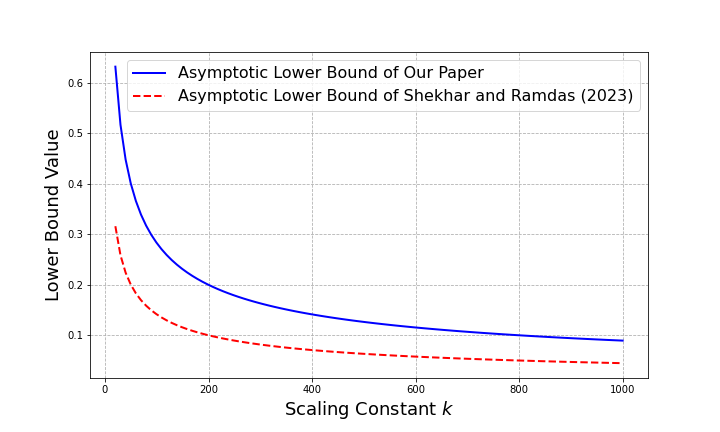}
    \caption{Comparison of our asymptotic lower bound given in Theorem \ref{thm_spef_lower_bound_width} as a function of $k$, with the lower bound presented in Proposition 4.3 of \cite{shekhar2023near} when $\lim_{\delta \to 0} \frac{N_\delta}{\log(1/\delta)} = k$. We assume that $\nu = N(0,1)$ with known variance.}
    \label{fig_lower}
\end{figure}

Observe that Proposition 4.3 of \cite{shekhar2023near} presents a non-asymptotic lower bound on the width of any CI for a given $N$ and $\delta$. We first compare our asymptotic lower bound with the one presented by \cite{shekhar2023near} in the asymptotic regime where $\delta \to 0$ and $N_\delta \to \infty$. In this regime, we plot both lower bounds versus the scaling constant $k$, where $\lim_{\delta \to 0} \tfrac{N_\delta}{\log(1/\delta)} = k$, for the case when $\nu = N(0,1)$ with known variance. The results demonstrate that our asymptotic lower bound is twice as high as that of \cite{shekhar2023near}. This is shown as a function of the scaling constant in Figure~\ref{fig_lower}.

The explanation for the result observed in Figure~\ref{fig_lower} is as follows. For the Gaussian case with known variance $\sigma$, we have $\mu^*_R(\mu) = \mu + \sigma\sqrt{2/k}$ and $\mu^*_L(\mu) = \mu - \sigma\sqrt{2/k}$. Thus, our lower bound is $2\sigma\sqrt{2/k}$, while the bound from Shekhar and Ramdas (2023) is $\sigma\sqrt{2/k}$.  

Hence, the factor of $2$ applies in the Gaussian case. For other distributions, this ratio may vary; however, our bound remains asymptotically tighter. Moreover, recall that we have proven that the policy $\pi_1$ achieves the limiting width in the sufficient learning regime, which shows that our lower bound cannot be improved in the regime when $\delta \to 0$ and $N_\delta \to \infty$.

We now compare the performance of the method $\pi_1$ with three existing methods: the Hoeffding-based CI, the Empirical Bernstein (EB) CI (see \cite{maurer2009empirical}), and the betting-based hedged CI (see \cite{waudby2024estimating}). Experiments were conducted on Bernoulli distributions with means of 0.6 and 0.9 and $\delta$ is set to be $1 \%$. For the betting method, we vary the discretization parameter \( m \in \{1000, 3000, 5000\} \) (needed to discretize the [0,1] space). For each configuration, we generated 1000 i.i.d. datasets of size \( N \in \{2000, 3000\} \), computed CI widths, and report the average width (max 95$\%$ CI width: $0.0001$). Our method consistently yields narrower CIs than both Hoeffding and EB across settings. While betting-based hedged CIs slightly outperform ours at \(N = 2000\), performance equalizes by \(N = 3000\). Increasing discretization \(m\) offers diminishing returns beyond \(m = 5000\). 
Further, our method remains computationally more efficient as compared to the betting-based hedged CIs due to its ease of computation.

\begin{table}[ht]
\centering
\caption{Average CI width for Bernoulli distributions with means 0.6 and 0.9 for $\delta = 0.01$}
\resizebox{\columnwidth}{!}{%
\begin{tabular}{ccccccc}
\toprule
\(N\) & $\pi_1$ & Betting ($m{=}1000$) & Betting ($m{=}3000$) & Betting ($m{=}5000$) & Hoeffding & EB \\
\midrule
\multicolumn{7}{c}{\textit{Mean = 0.6}} \\
2000 & 0.0712 & 0.0603 & 0.0596 & 0.0595 & 0.0728 & 0.0898 \\
3000 & 0.0582 & 0.0592 & 0.0585 & 0.0584 & 0.0594 & 0.0712 \\
\midrule
\multicolumn{7}{c}{\textit{Mean = 0.9}} \\
2000 & 0.0436 & 0.0378 & 0.0371 & 0.0369 & 0.0728 & 0.0606 \\
3000 & 0.0356 & 0.0370 & 0.0363 & 0.0361 & 0.0594 & 0.0473 \\
\bottomrule
\end{tabular}%
}
\end{table}

\section{CONCLUSION AND FUTURE DIRECTIONS} \label{sec_conclusion}
 In this paper, we explored the construction of confidence intervals (CIs) for the mean of a probability distribution, focusing on their minimum achievable limiting width in three distinct learning regimes. By characterizing these regimes, we provided a comprehensive understanding of the inherent limitations and potential of CI construction methods in both parametric and non-parametric settings. Finally, we demonstrated the asymptotic optimality of KL-based CI construction methods and extended our results to practically relevant scenarios where sample collection incurs random costs. While non-asymptotic results can capture more “granularity” and reveal finer structural properties, nonetheless, our asymptotic analysis offers more concise and elegant insights that have independent value. In particular, our stable policies assumption is asymptotic, is true quite generally and greatly simplifies the analysis. Future research directions could include extending the framework to other statistical estimation problems, such as variance and quantile estimation. At last, extending our results to higher-dimensional parameters, such as the mean of a multivariate distribution, likewise represents a compelling and nontrivial avenue for future research.

\bibliography{references}

\clearpage
\appendix
\thispagestyle{empty}

\onecolumn
\aistatstitle{Supplementary Material}

\noindent \textbf{Discussion on stable policies:} 
Fix a distribution $\nu$ with mean $\mu$. Recall a policy $\pi$ is called stable if    
$$
\lim_{\delta \to 0 } \widehat{\mu}^{\pi}_L(N_\delta,\delta) \overset{p}{\to} \mu_{L}^{\pi}(\nu), \quad
\lim_{\delta \to 0 } \widehat{\mu}^{\pi}_R(N_\delta,\delta) \overset{p}{\to} \mu_{R}^{\pi}(\nu),
$$
where $\mu_{L}^{\pi}(\nu) \leq  \mu$ and   $\mu_{R}^{\pi}(\nu) \geq  \mu$ are constants ($-\infty$ and $\infty$ included).
For instance, a policy using a symmetric confidence interval based on the central limit theorem inherently meets this assumption. Furthermore, for bounded outcome distributions, classical confidence interval approaches based on Hoeffding's and Bernstein's inequalities, as well as Maurer and Pontil's Empirical Bernstein (MP-EB) CI \cite{maurer2009empirical}, also satisfy the stability assumption. To verify this, we analyze the asymptotic behavior of the CI boundaries as $N_\delta \to \infty$ and $\delta \to 0$ for three CI constructions. The key observation is that the CI boundaries converge almost surely to constants depending on the relationship between $N_\delta$ and $\log(1/\delta)$. Specifically, we consider three cases based on $\lim_{\delta \to 0} \frac{N_\delta}{\log(1/\delta) }= k$, where $k \in \{0, (0,\infty), \infty\}$.

\textbf{Hoeffding CI.}
$
C_N^{(H)} = \left[\widehat{\mu}_N \pm \frac{w_N^{(H)}}{2}\right]
:= \left[\widehat{\mu}_N - \frac{w_N^{(H)}}{2}, \widehat{\mu}_N + \frac{w_N^{(H)}}{2}\right],
$
where $\widehat{\mu}_N = \frac{\sum_{i=1}^{N}X_i}{N}$ and $w_N^{(H)} = 2\sqrt{\frac{\log(2/\delta)}{2N}}$.
For Hoeffding's CI, by the strong law of large numbers, $\widehat{\mu}_N \overset{a.s.}{\to} \mu$ as $N \to \infty$. The width term satisfies:
\begin{itemize}
\item If $k = \infty$: $w_N^{(H)} = 2\sqrt{\frac{\log(2/\delta)}{2N}} \to 0$, thus $\mu_L^{\pi}(\nu) = \mu_R^{\pi}(\nu) = \mu$.
\item If $k \in (0,\infty)$: $w_N^{(H)} \to \sqrt{\frac{2}{k}}$ (constant), thus $\mu_L^{\pi}(\nu) = \mu - \frac{1}{2}\sqrt{\frac{2}{k}}$ and $\mu_R^{\pi}(\nu) = \mu + \frac{1}{2}\sqrt{\frac{2}{k}}$.
\item If $k = 0$: $w_N^{(H)} \to \infty$, thus $\mu_L^{\pi}(\nu) = -\infty$ and $\mu_R^{\pi}(\nu) = +\infty$.
\end{itemize}

\textbf{Bernstein CI.}
$
C_N^{(B)} = \left[\widehat{\mu}_N \pm \frac{w_N^{(B)}}{2}\right],
\quad \text{where } w_N^{(B)} = 2\sigma\sqrt{\frac{2\log(2/\delta)}{N}} + \frac{4\log(2/\delta)}{3N}.
$
For Bernstein's CI (see \cite{shekhar2023near} for above formulation), the analysis is similar. 
\begin{itemize}
\item If $k = \infty$: $w_N^{(B)} \to 0$, thus $\mu_L^{\pi}(\nu) = \mu_R^{\pi}(\nu) = \mu$.
\item If $k \in (0,\infty)$: $w_N^{(B)} \to 2\sigma\sqrt{\frac{2}{k}} + \frac{4}{3k}$ (constant), thus $\mu_L^{\pi}(\nu) = \mu - \sigma\sqrt{\frac{2}{k}} - \frac{2}{3k}$ and $\mu_R^{\pi}(\nu) = \mu + \sigma\sqrt{\frac{2}{k}} + \frac{2}{3k}$.
\item If $k = 0$: $w_N^{(B)} \to \infty$, thus $\mu_L^{\pi}(\nu) = -\infty$ and $\mu_R^{\pi}(\nu) = +\infty$.
\end{itemize}

\textbf{Maurer $\&$ Pontil's Empirical Bernstein (MP-EB) CI.}
$
C_N^{\text{(MP-EB)}} = \left[\widehat{\mu}_N \pm \frac{w_N^{\text{(MP-EB)}}}{2}\right],
$
where
$
w_N^{\text{(MP-EB)}} = 2\widehat{\sigma}_N\sqrt{\frac{2\log(4/\delta)}{N}} + \frac{14\log(4/\delta)}{3(N-1)},
\quad \text{and } \widehat{\sigma}^2{_N} = \frac{\sum_{i=1}^{N}(X_i - \widehat{\mu}_N)^2}{N-1}.
$
For the MP-EB CI, note that $\widehat{\sigma}_{N_\delta} \overset{a.s.}{\to} \sigma$ as $N_\delta \to \infty$. Hence,
\begin{itemize}
\item If $k = \infty$: $w_N^{\text{(MP-EB)}} \to 0$, thus $\mu_L^{\pi}(\nu) = \mu_R^{\pi}(\nu) = \mu$.
\item If $k \in (0,\infty)$: $w_N^{\text{(MP-EB)}} \to 2\sigma\sqrt{\frac{2}{k}} + \frac{14}{3k}$ (constant), thus $\mu_L^{\pi}(\nu) = \mu - \sigma\sqrt{\frac{2}{k}} - \frac{7}{3k}$ and $\mu_R^{\pi}(\nu) = \mu + \sigma\sqrt{\frac{2}{k}} + \frac{7}{3k}$.
\item If $k =0 $: $w_N^{\text{(MP-EB)}} \to \infty$, thus $\mu_L^{\pi}(\nu) = -\infty$ and $\mu_R^{\pi}(\nu) = +\infty$.
\end{itemize}

In all these cases, the CI boundaries converge almost surely to constants as $N_\delta \to \infty$ and $\delta \to 0$, confirming that these CI-based policies satisfy the stability assumption.

\section{Proofs of results in Section \ref{sec_lower_bound}.}
\label{appen_lower_bound}
    
\noindent \textbf{Proof of Theorem \ref{thm_spef_lower_bound_width}.}
Consider any $\pi \in \Pi_{\rm CI}^{s}$. As the policy $\pi$ is stable, it follows that for a given distribution $\nu \in \mathbf{S}$ with mean $\mu$, we have,
\begin{equation} \label{appen_0}  
\lim_{\delta \to 0} \hat{\mu}_L^{\pi}(N_\delta,\delta) \overset{p}{\to}  \mu_{L}^{\pi}(\mu) \textrm{ and } \lim_{\delta \to 0} \hat{\mu}_R^{\pi}(N_\delta,\delta) \overset{p}{\to}  \mu_{R}^{\pi}(\mu) .
\end{equation} 
 We now define the set of alternate environments $K(\mu_L^\pi(\mu),\mu_R^\pi(\mu)) = K_1(\mu_L^\pi(\mu)) \cup K_2(\mu_R^\pi(\mu))$, $  K_1(\mu_L^\pi(\mu)) = \{ \tilde{\nu}: \tilde{\nu} \in \mathbf{S}, \tilde{\mu} < \mu_L^\pi(\mu)\}$ and  $K_2(\mu_R^\pi(\mu)) = \{ \tilde{\nu}: \tilde{\nu} \in \mathbf{S}, \tilde{\mu} > \mu_R^\pi(\mu)\}$ . 

 Using the data processing inequality for a given $\delta \in (0,1)$ and any alternate environment $\tilde{\nu}$ with mean $\tilde{\mu}$, we have
	\begin{equation} \label{appen_1}
N_\delta \cdot d(\mu, \tilde{\mu})  \geq \sup_{\mathcal{E}_\delta \in \mathcal{F}_{N_\delta}} \phi(\mathbb{P}_{\nu}(\mathcal{E}_\delta),
	\mathbb{P}_{\tilde{\nu}}(\mathcal{E}_\delta)),
  	\end{equation}
where  $\phi(p_1,p_2) \triangleq  p_1 \log\frac{p_1}{p_2} +(1-p_1)\log \left(\frac{1-p
_1}{1-p_2}\right)$ for $p_1, p_2 \in (0,1)$. One can use Lemma 0.1 in \cite{kaufmann2020contributions} to get \eqref{appen_1} from the data processing inequality.

We first utilize \eqref{appen_1} for $\tilde{\nu} \in K_1(\mu_L^\pi(\mu))$ and $\mathcal{E}_\delta = \{\tilde{\mu} \notin [\widehat{\mu}^{\pi}_L(N_\delta,\delta), \widehat{\mu}^{\pi}_R(N_\delta,\delta)]\}$. A similar approach would follow for $\tilde{\nu} \in K_2(\mu_R^\pi(\mu))$.

Since $\pi \in \Pi_{\rm CI}^{s}$, using the definition of CI, we get, $\mathbb{P}_{\tilde{\nu}}(\mathcal{E}_\delta) \leq \delta$. Now observe that,
$$\mathbb{P}_{\nu}(\mathcal{E}_\delta) \geq \mathbb{P}_{\nu} \{\tilde{\mu} < \widehat{\mu}^{\pi}_L(N_\delta,\delta)\} . $$
Taking the limit of $\delta \to 0$ on both sides and using \eqref{appen_0}, we get,

 $$ \lim_{\delta \to 0} \mathbb{P}_{\nu}(\mathcal{E}_\delta) =1 .$$

Hence using the definition of $\phi(\cdot)$, we get,
 
\begin{equation} \label{appen_2} 
\liminf_{\delta \to 0} \frac{ \phi(\mathbb{P}_{\nu}(\mathcal{E}_\delta),  \mathbb{P}_{\tilde{\nu}}(\mathcal{E}_\delta))}{\log (1/ \delta)} \geq 1,
\end{equation}

Hence using \eqref{appen_1} and \eqref{appen_2} it follows that for all $\tilde{\nu} \in K_1(\mu_L^\pi(\mu))$ with mean $\tilde{\mu}$, we have $ \liminf_{\delta \to 0 }\frac{N_\delta}{\log(1/\delta)} \cdot d(\mu, \tilde{\mu})  \geq 1.$ Using a similar analysis for  $\tilde{\nu} \in K_2(\mu_R^\pi(\mu))$ with mean $\tilde{\mu}$, we have a similar result, i.e., $ \liminf_{\delta \to 0 }\frac{N_\delta}{\log(1/\delta)} \cdot d(\mu, \tilde{\mu})  \geq 1.$ Hence combining both, we get that for all $\tilde{\nu} \in K(\mu_L^\pi(\mu), \mu_R^\pi(\mu))$ with mean $\tilde{\mu}$ following holds,
$$ \liminf_{\delta \to 0 }\frac{N_\delta}{\log(1/\delta)}  d(\mu, \tilde{\mu})  \geq 1.$$

It follows that,
\begin{equation} \label{appen_3}
    \liminf_{\delta \to 0 }\frac{N_\delta}{\log(1/\delta)}  d(\mu, \mu_L^\pi(\mu) -\eta)  \geq 1 \textrm{ and } \liminf_{\delta \to 0 }\frac{N_\delta}{\log(1/\delta)}  d(\mu, \mu_R^\pi(\mu) +\eta)  \geq 1 ,
\end{equation}
where $\eta >0$ is any small positive number.
Now we consider the case when $\liminf_{\delta \to 0 }\frac{N_\delta}{\log(1/\delta)} = 0$. In this case, \eqref{appen_3} can not hold for any $\mu_L^\pi(\mu) \in \mathbf{O}$ or $\mu_R^\pi(\mu) \in \mathbf{O}$. Hence, we can define $\mu_L^\pi(\mu) = \underline{\mu}$ and  $\mu_R^\pi(\mu) = \overline{\mu}$ for any policy $\pi \in \Pi^{s}_{\rm CI}$. This completes the proof of part (a).

Now we consider the case when $\liminf_{\delta \to 0 }\frac{N_\delta}{\log(1/\delta)} = k
$ for $k \in (0,\infty)$. In this case using \eqref{appen_3} and taking $\eta \to 0$, we get, 
$$    d(\mu, \mu_L^\pi(\mu))  \geq \frac{1}{k} \textrm{ and }  d(\mu, \mu_R^\pi(\mu))  \geq \frac{1}{k} . 
$$
Using the strict quasi-convexity of $d(\mu,\cdot)$, we get the desired result. This completes the proof.

\hfill $\Box$

\section{Proofs of results in Section \ref{sec_kl}.}

\noindent \textbf{Proof of Theorem \ref{thm_kl}.}

 We first show that our policy $\pi_1 \in \Pi_{\rm CI}^{s}$. To prove this, we first show that,
\begin{equation} \label{appen_7}
     \forall n \in \mathbb{N} : \mathbb{P}_{\nu}(\mu \in [\widehat{\mu}^{\pi_1}_L(n,\delta), \widehat{\mu}^{\pi_1}_R(n,\delta)]) \geq 1-\delta.
\end{equation}

Observe that, for any $n \in \mathbb{N}$,  it suffices to show that

$$\mathbb{P}_{\nu}(\mu \notin [\widehat{\mu}^{\pi_1}_L(n,\delta), \widehat{\mu}^{\pi_1}_R(n,\delta)])  \leq \mathbb{P}_{\nu}(n d(\hat{\mu}_n,\mu) \geq \beta(\delta)) = \mathbb{P}_{\nu}(n d(\hat{\mu}_n,\mu) \geq  \log(2/\delta)) \leq \delta.$$

Using the result from \cite{garivier2011kl}, we know that for $\nu \in \mathbf{S}$,
\begin{equation} \label{appen_4}
    d(\hat{\mu}_n,\mu) = I(\hat{\mu}_n) =  \sup_{\lambda \in \mathbb{R}} \lambda \hat{\mu}_n - \log \mathbb{E}_\nu [e^{\lambda X}].
\end{equation}
Here $I(\cdot)$ is the good rate function used in large deviation.
In the above, $X$ is distributed according to $\nu$.

Observe that we need to show (this can be derived from Lemma 4 in \cite{menard2017minimax} or Theorem 4 in \cite{busa2019optimal} as well),
$$\mathbb{P}_{\nu}(n d(\hat{\mu}_n,\mu) \geq  \log(2/\delta)) \leq \delta.$$
It suffices to show that,  
$$\mathbb{P}_{\nu}(n d(\hat{\mu}_n,\mu) \geq  \log(2/\delta) \cap  \hat{\mu}_n < \mu) \leq \frac{\delta }{2} \textrm{ and } \mathbb{P}_{\nu}(n d(\hat{\mu}_n,\mu) \geq  \log(2/\delta) \cap  \hat{\mu}_n \geq \mu) \leq \frac{\delta }{2}.$$

Now we show that
\begin{equation} \label{appen_5}
  \mathbb{P}_{\nu}(n d(\hat{\mu}_n,\mu) \geq  \log(2/\delta) \cap  \hat{\mu}_n < \mu) \leq \frac{\delta }{2}.  
\end{equation}
Let $z \in (\underline{\mu},\mu)$ such that, $I(z) = d(z,\mu) = \frac{\log(2/\delta)}{n}.$ Using the property of good rate function we know that, for $z< \mu$, there exists a $\lambda(z) <0$ such that the following holds
$$I(z)= d(z,\mu) = \lambda(z)z -\log \mathbb{E}_\nu [e^{\lambda(z) X}]. $$

Now under the event $\{n d(\hat{\mu}_n,\mu) \geq  \log(2/\delta) \cap  \hat{\mu}_n < \mu\}$, we have,

$$\hat{\mu}_n \leq z.$$ Hence using the fact that $\lambda(z) <0$, we get,

$$\{n d(\hat{\mu}_n,\mu) \geq  \log(2/\delta) \cap  \hat{\mu}_n < \mu\} \implies  \left\{\lambda(z) \hat{\mu}_n - \log \mathbb{E}_\nu[e^{\lambda(z) X}] \geq  d(z,\mu) = \frac{\log(2/\delta)}{n} \right\}.$$

Observe that, $\hat{\mu}_n = \frac{\sum_{i=1}^{n} X_i}{n}$, to prove \eqref{appen_5}, it suffices to show that,
$$\mathbb{P}_{\nu}\left(\lambda(z) \sum_{i=1}^{n}X_i - n \log \mathbb{E}_\nu[e^{\lambda(z) X}]\geq  \log(2/\delta) \right)\leq \frac{\delta }{2}.$$

Observe that  $\mathbb{E}_{\nu}[e^{\lambda(z) \sum_{i=1}^{n}X_i - n \log \mathbb{E}_\nu[e^{\lambda(z) X}]}] = 1$, hence using Markov's inequality, we get the desired result.

Similarly, one can show that, 

$$\mathbb{P}_{\nu}(n d(\hat{\mu}_n,\mu) \geq  \log(2/\delta) \cap  \hat{\mu}_n \geq \mu) \leq \frac{\delta }{2}.$$

This completes the proof of \eqref{appen_7}.

 Observe that from the definitions of $\widehat{\mu}^{\pi_1}_L(n,\delta)$, $\widehat{\mu}^{\pi_1}_R(n,\delta)$ and the strict quasi-convexity of $d(\mu,\cdot)$, we get,

$$d(\hat{\mu}_{N_\delta}, \widehat{\mu}^{\pi_1}_L(N_\delta,\delta)) = d(\hat{\mu}_{N_\delta}, \widehat{\mu}^{\pi_1}_R(N_\delta,\delta)) = \frac{\log(2/\delta)}{N_\delta}.$$

Recall $N_\delta \to \infty$, $ \lim_{\delta \to 0} \hat{\mu}_{N_\delta} \to \mu$ almost surely from strong law of large numbers.
Taking $\delta \to 0$ and using the joint-continuity of $d(\mu,x)$ in $(\mu,x)$, we get that $\pi_1 \in \Pi_{\rm CI}^{s}$ and part (b), (c) of this theorem holds. This completes the proof.

\hfill $\Box$

\section{Results and proofs for the rate of convergence analysis in complete learning regime.}\label{sub_sec_rate_analysis}

We first introduce a set of policies, denoted by $\Pi_{\rm CI}^{sr} \subseteq \Pi_{\rm CI}^{s}$ for which the analysis is valid. For $\pi \in \Pi_{\rm CI}^{sr}$ and a given $\nu \in \mathbf{S}$ with mean $\mu$,  we assume that the following holds as $ \lim_{\delta \to 0} \frac{N_\delta}{\log(1/\delta)} = \infty$:

 $$ \lim_{\delta \to 0}\frac{\mu- \widehat{\mu}^{\pi}_L(N_\delta,\delta) }{\sqrt{\frac{\log(1/\delta)}{N_\delta}}} \overset{p}{\to} \theta_L^{\pi}(\mu) \textrm{ and } \lim_{\delta \to 0}\frac{ \widehat{\mu}^{\pi}_R(N_\delta,\delta) - \mu }{\sqrt{\frac{\log(1/\delta)}{N_\delta}}} \overset{p}{\to} \theta_R^{\pi}(\mu),$$
 where, $\theta_L^{\pi}(\mu)$ and $\theta_R^{\pi}(\mu)$ are the non-negative constants.

It is worth noting that the above assumption can be interpreted as rate stability of policies under the complete learning regime. The example of stable policies discussed at the start of the supplementary material also satisfies this rate stability assumption.

We now state a lower bound on the rate of convergence of the width to zero.

\begin{theorem} \label{thm_complete_learning}
 Fix a $\nu \in \mathbf{S}$ with mean $\mu$. If $ \lim_{\delta \to 0} \frac{N_\delta}{\log(1/\delta)} = \infty$, then for any $\pi \in \Pi^{sr}_{\rm CI}$, the following holds:
  $$  \lim_{\delta \to 0 }\frac{\widehat{\mu}^{\pi}_R(N_\delta,\delta) - \widehat{\mu}^{\pi}_L(N_\delta,\delta)}{\sqrt{\frac{\log(1/\delta)}{N_\delta}}}   \overset{p}{\to}q(\mu)\geq   \sqrt{8} \cdot \sigma(\mu),$$
where $\sigma(\mu)$ is the standard deviation of the $\nu$. Further, $q(\mu)$ is some non-negative function. 
\end{theorem}

\noindent\textbf{Proof of Theorem \ref{thm_complete_learning}}.
Recall that from \eqref{appen_1}, for a given $\delta \in (0,1)$ and any alternate environment $\tilde{\nu}$ with mean $\tilde{\mu}$, we have
	\begin{equation} 
N_\delta \cdot d(\mu, \tilde{\mu})  \geq \sup_{\mathcal{E}_\delta \in \mathcal{F}_{N_\delta}} \phi(\mathbb{P}_{\nu}(\mathcal{E}_\delta),
	\mathbb{P}_{\tilde{\nu}}(\mathcal{E}_\delta)). \nonumber 
  	\end{equation}

Choose $\tilde{\mu} = \mu - (\theta_L^{\pi}(\mu) + \eta) \sqrt{\frac{\log(1/\delta)}{N_\delta}}$. Since we know that $d(\mu,x)$ is twice continuously differentiable in $x$ and $\frac{\partial d(\mu,x)}{\partial x}|_{x = \mu} = 0$. Hence, using the Taylor series expansion of $d(\mu, \cdot)$ in the above equation, we get
$$N_\delta \frac{I(c_1) (\theta_L^{\pi}(\mu) + \eta)^2 \log(1/\delta)}{2N_\delta} \geq \sup_{\mathcal{E}_\delta \in \mathcal{F}_{N_\delta}} \phi(\mathbb{P}_{\nu}(\mathcal{E}_\delta),
\mathbb{P}_{\tilde{\nu}}(\mathcal{E}_\delta)).$$

Here $I(c) = \frac{\partial^2 d(\mu,x)}{\partial x^2}|_{x = c}$ and $c_1 \in (\tilde{\mu}, \mu)$ is a constant. Now we choose $\mathcal{E}_\delta = \{\tilde{\mu} \notin [\widehat{\mu}^{\pi}_L(N_\delta,\delta), \widehat{\mu}^{\pi}_R(N_\delta,\delta)]\}$ and take $\delta \to 0$, we get

\begin{equation} \label{appendix_11}
   \frac{ I(\mu) (\theta_L^{\pi}(\mu) + \eta)^2 }{2}\geq \lim_{\delta \to 0} \frac{\phi(\mathbb{P}_{\nu}(\mathcal{E}_\delta),
	\mathbb{P}_{\tilde{\nu}}(\mathcal{E}_\delta))}{\log(1/\delta)}.
\end{equation}
 Now consider the following:
$$\mathbb{P}_{\nu}(\mathcal{E}_\delta) = \mathbb{P}_{\nu} (\tilde{\mu} \notin [\widehat{\mu}^{\pi}_L(N_\delta,\delta), \widehat{\mu}^{\pi}_R(N_\delta,\delta)]) .$$
It follows that,
$$\mathbb{P}_{\nu}(\mathcal{E}_\delta) \geq  \mathbb{P}_{\nu} (\tilde{\mu} < \widehat{\mu}^{\pi}_L(N_\delta,\delta)) .$$

Substituting the definition of $\tilde{\mu}$, we get,

$$\mathbb{P}_{\nu}(\mathcal{E}_\delta) \geq \mathbb{P}_{\nu} \left(\mu -\widehat{\mu}^{\pi}_L(N_\delta,\delta) < (\theta_L^{\pi}(\mu) + \eta) \sqrt{\frac{\log(1/\delta)}{N_\delta}} \right) .$$

It can be re-written as,

$$ \lim_{\delta \to 0}\mathbb{P}_{\nu}(\mathcal{E}_\delta) \geq \lim_{\delta \to 0} \mathbb{P}_{\nu} \left(  \frac{\mu -\widehat{\mu}^{\pi}_L(N_\delta,\delta)}{\sqrt{\frac{\log(1/\delta)}{N_\delta}}} <  (\theta_L^{\pi}(\mu) + \eta) \right) .$$

Now using the fact that $\pi \in \Pi_{\rm CI}^{sr}$, we know that,  $\lim_{\delta \to 0} \frac{ \mu  -\widehat{\mu}^{\pi}_L(N_\delta,\delta)}{\sqrt{\frac{\log(1/\delta)}{N_\delta}}} \overset{p}{\to } \theta_{L}^{\pi}(\mu)$.  Hence it follows that, for any $\eta> 0$, we have,

$$ \lim_{\delta \to 0}\mathbb{P}_{\nu}(\mathcal{E}_\delta) \geq \lim_{\delta \to 0} \mathbb{P}_{\nu} \left(  \frac{\mu -\widehat{\mu}^{\pi}_L(N_\delta,\delta)}{\sqrt{\frac{\log(1/\delta)}{N_\delta}}} <  (\theta_L^{\pi}(\mu) + \eta) \right) = 1 .$$

Observe that we get trivially, $\mathbb{P}_{\tilde{\nu}}(\mathcal{E}_\delta) \leq \delta.$ Substituting the value of $\mathbb{P}_{\tilde{\nu}}(\mathcal{E}_\delta)$ and $ \lim_{\delta \to 0}\mathbb{P}_{\nu}(\mathcal{E}_\delta)$ in \eqref{appendix_11}, we get that,

$$  \frac{ I(\mu) (\theta_L^{\pi}(\mu) + \eta)^2 }{2}\geq \lim_{\delta \to 0} \frac{\phi(\mathbb{P}_{\nu}(\mathcal{E}_\delta),
	\mathbb{P}_{\tilde{\nu}}(\mathcal{E}_\delta))}{\log(1/\delta)} \geq 1.$$
   Taking $\eta \to 0$ and the fact that $I(\mu) = \frac{1}{\sigma^2(\mu)}$, we get, 

   $$ \theta_L^{\pi}(\mu) \geq  \sqrt{2} \cdot \sigma(\mu).$$

   Similarly, we get, 
      $$ \theta_R^{\pi}(\mu) \geq   \sqrt{2} \cdot \sigma(\mu).$$

      Now observe that, 
      $$\lim_{\delta \to 0 }\frac{\widehat{\mu}^{\pi}_R(N_\delta,\delta) - \widehat{\mu}^{\pi}_L(N_\delta,\delta)}{\sqrt{\frac{\log(1/\delta)}{N_\delta}}} =  \lim_{\delta \to 0 } \frac{\widehat{\mu}^{\pi}_R(N_\delta,\delta) - \mu}{\sqrt{\frac{\log(1/\delta)}{N_\delta}}} + \lim_{\delta \to 0 } \frac{\mu - \widehat{\mu}^{\pi}_L(N_\delta,\delta)}{\sqrt{\frac{\log(1/\delta)}{N_\delta}}}\overset{p}{\to} \theta_L^{\pi}(\mu) + \theta_R^{\pi}(\mu) .  $$
      Using the fact shown above that, $\theta_L^{\pi}(\mu) \geq  \sqrt{2} \cdot \sigma(\mu)$ and $\theta_R^{\pi}(\mu) \geq  \sqrt{2} \cdot \sigma(\mu)$, we get the desired result.

\hfill $\Box$

Now we state the result which shows that our policy $\pi_1$ has the fastest rate of convergence of width to zero in the complete learning regime.

\begin{theorem}  \label{last}
  The policy $\pi_1$ has the following properties:
    
    a) $\pi_1 \in \Pi^{sr}_{\rm CI}$.
    b) If $\lim_{\delta \to 0} \frac {N_\delta} {\log(1/\delta)}= \infty$, then we have, 
 $$   \lim_{\delta \to 0 }\frac{\widehat{\mu}^{\pi_1}_R(N_\delta,\delta) - \widehat{\mu}^{\pi_1}_L(N_\delta,\delta)}{\sqrt{\frac{\log(1/\delta)}{N_\delta}}}   \overset{p}{\to} \sqrt{8} \cdot \sigma(\mu).$$

\end{theorem}

\noindent \textbf{Proof of the Theorem \ref{last}.}

Observe that from the definition of $\widehat{\mu}^{\pi_1}_L(n,\delta)$, $\widehat{\mu}^{\pi_1}_R(n,\delta)$ and the strict quasi-convexity of $d(\mu,\cdot)$, we get,
\begin{equation} \label{appendix_12}
    d(\hat{\mu}^{\pi_1}_{N_\delta}, \widehat{\mu}^{\pi_1}_L(N_\delta,\delta)) = d(\hat{\mu}_{N_\delta}^{\pi_1}, \widehat{\mu}^{\pi_1}_R(N_\delta,\delta)) = \frac{\log(2/\delta)}{N_\delta}.
\end{equation}

First, note that $\hat{\mu}^{\pi_1}_{N_\delta}$ is the sample average based on $N_\delta$ observations. Since $N_\delta \to \infty$ as $\delta \to 0$, the central limit theorem applies. In particular,
\[
\sqrt{N_\delta}\,(\hat{\mu}^{\pi_1}_{N_\delta} - \mu)
\;\xrightarrow{d}\;
\mathcal{N}(0, \sigma^2(\mu)),
\]

This implies that
\[
\sqrt{N_\delta}\,(\hat{\mu}^{\pi_1}_{N_\delta} - \mu) = O_p(1).
\]

Now consider
\[
\frac{\sqrt{N_\delta} (\mu-\hat{\mu}^{\pi_1}_{N_\delta})}{ \sqrt{\log(1/\delta)}} 
=
\left(\sqrt{N_\delta} (\mu-\hat{\mu}^{\pi_1}_{N_\delta})\right)
\cdot
\frac{1}{\sqrt{\log(1/\delta)}}.
\]

Since $\delta \to 0$, we have $\log(1/\delta) \to \infty$. Therefore, using Slutsky's theorem (or equivalently, the fact that $O_p(1)\cdot o(1) \xrightarrow{p} 0$), we obtain

\begin{equation} \label{new_chat}
    \lim_{\delta \to 0} \frac{\sqrt{N_\delta} (\mu-\hat{\mu}^{\pi_1}_{N_\delta})}{ \sqrt{\log(1/\delta)}} \overset{p}{\to} 0.
\end{equation}

Using Taylor's theorem,  there exist random points
\[
c_{\delta,L}\in 
\bigl(\widehat{\mu}^{\pi_1}_L(N_\delta,\delta),
\hat{\mu}^{\pi_1}_{N_\delta}\bigr),
\qquad
c_{\delta,R}\in 
\bigl(\hat{\mu}^{\pi_1}_{N_\delta},
\widehat{\mu}^{\pi_1}_R(N_\delta,\delta)\bigr)
\]
such that
\[
d(\hat{\mu}^{\pi_1}_{N_\delta},
\widehat{\mu}^{\pi_1}_L(N_\delta,\delta))
=
\frac{1}{2}
\left.
\frac{\partial^2}{\partial x^2}
d(\hat{\mu}^{\pi_1}_{N_\delta},x)
\right|_{x=c_{\delta,L}}
\left(\hat{\mu}^{\pi_1}_{N_\delta}
-
\widehat{\mu}^{\pi_1}_L(N_\delta,\delta)\right)^2,
\]
and
\[
d(\hat{\mu}^{\pi_1}_{N_\delta},
\widehat{\mu}^{\pi_1}_R(N_\delta,\delta))
=
\frac{1}{2}
\left.
\frac{\partial^2}{\partial x^2}
d(\hat{\mu}^{\pi_1}_{N_\delta},x)
\right|_{x=c_{\delta,R}}
\left(\widehat{\mu}^{\pi_1}_R(N_\delta,\delta)
-
\hat{\mu}^{\pi_1}_{N_\delta}\right)^2.
\]
Moreover, in the complete learning regime, $
\frac{\log(2/\delta)}{N_\delta}\to 0.$
Therefore, using \eqref{appendix_12}, both 
$\widehat{\mu}^{\pi_1}_L(N_\delta,\delta)$ and 
$\widehat{\mu}^{\pi_1}_R(N_\delta,\delta)$ converge in probability to 
$\hat{\mu}^{\pi_1}_{N_\delta}$, and since 
$\hat{\mu}^{\pi_1}_{N_\delta}\overset{p}{\to}\mu$, we also have
\[
c_{\delta,L}\overset{p}{\to}\mu,
\qquad
c_{\delta,R}\overset{p}{\to}\mu.
\]

By the continuity of $(\mu,x)\mapsto \frac{\partial^2}{\partial x^2}d(\mu,x)$ in a neighborhood of $(\mu,\mu)$, together with the continuous mapping theorem, we get
\[
\left.
\frac{\partial^2}{\partial x^2}
d(\hat{\mu}^{\pi_1}_{N_\delta},x)
\right|_{x=c_{\delta,L}}
\overset{p}{\to}
\frac{1}{\sigma^2(\mu)},
\qquad
\left.
\frac{\partial^2}{\partial x^2}
d(\hat{\mu}^{\pi_1}_{N_\delta},x)
\right|_{x=c_{\delta,R}}
\overset{p}{\to}
\frac{1}{\sigma^2(\mu)}.
\]
Combining this with \eqref{appendix_12}, we obtain

$$ \lim_{\delta \to 0}\frac{\hat{\mu}^{\pi_1}_{N_\delta}- \widehat{\mu}^{\pi_1}_L(N_\delta,\delta)}{\sqrt{\frac{\log(1/\delta)}{N_\delta}}} \overset{p}{\to} \sqrt{2} \cdot \sigma(\mu) \textrm{ and } \lim_{\delta \to 0}\frac{ \widehat{\mu}^{\pi_1}_R(N_\delta,\delta) - \hat{\mu}^{\pi_1}_{N_\delta} }{\sqrt{\frac{\log(1/\delta)}{N_\delta}}} \overset{p}{\to} \sqrt{2} \cdot \sigma(\mu).$$

 Using \eqref{new_chat} and the above, we get,

\begin{equation} \label{chat_new_2}
    \lim_{\delta \to 0}\frac{\mu- \widehat{\mu}^{\pi_1}_L(N_\delta,\delta)}{\sqrt{\frac{\log(1/\delta)}{N_\delta}}} \overset{p}{\to} \sqrt{2} \cdot \sigma(\mu) \textrm{ and } \lim_{\delta \to 0}\frac{ \widehat{\mu}^{\pi_1}_R(N_\delta,\delta) - \mu }{\sqrt{\frac{\log(1/\delta)}{N_\delta}}} \overset{p}{\to} \sqrt{2} \cdot \sigma(\mu).
\end{equation}
Using Theorem \ref{thm_kl}, we know that $\pi_1 \in \Pi_{\rm CI}^{s}$. Hence, using \eqref{chat_new_2}, we get that $\pi_1 \in \Pi_{\rm CI}^{sr}$. Further, it follows that, 

$$ \lim_{\delta \to 0 }\frac{\widehat{\mu}^{\pi_1}_R(N_\delta,\delta) - \widehat{\mu}^{\pi_1}_L(N_\delta,\delta)}{\sqrt{\frac{\log(1/\delta)}{N_\delta}}}  =  \lim_{\delta \to 0 }\frac{\widehat{\mu}^{\pi_1}_R(N_\delta,\delta) - \mu}{\sqrt{\frac{\log(1/\delta)}{N_\delta}}} + \lim_{\delta \to 0 }\frac{\mu - \widehat{\mu}^{\pi_1}_L(N_\delta,\delta)}{\sqrt{\frac{\log(1/\delta)}{N_\delta}}} \overset{p}{\to} \sqrt{8} \cdot \sigma(\mu).$$

This completes the proof.

\hfill $\Box$

\section{Proofs of results in Section \ref{sec_random}.}

\noindent \textbf{Proof of Theorem \ref{thm_spef_lower_bound_width_cost}.}

Recall that,
\begin{equation} 
    \tau_\delta = \sup \{ n \in \mathbb{Z}^{+}: \sum_{i = 1}^{n}c_i \leq C_\delta \}. \nonumber
\end{equation}

Using elementary renewal theorem for the renewal process, as $C_\delta \to \infty$, we get,

\begin{equation} \label{appe_10}
    \lim_{\delta \to 0}\frac{\mathbb{E}[\tau_\delta]}{C_\delta} = \frac{1}{\overline{c}} \textrm{ almost surely}.
\end{equation}

Consider any $\pi \in \hat{\Pi}_{\rm CI}^{s}$. As the policy $\pi$ is stable, it follows that for a given distribution $\nu \in \mathbf{S}$ with mean $\mu$ and a cost distribution $\mathcal{C}$ with mean $\overline{c}$, we have,
\begin{equation} \nonumber
\lim_{\delta \to 0} \hat{\mu}_L^{\pi}(\tau_\delta,\delta) \overset{p}{\to}  \mu_{L}^{\pi}(\mu) \textrm{ and } \lim_{\delta \to 0} \hat{\mu}_R^{\pi}(\tau_\delta,\delta) \overset{p}{\to}  \mu_{R}^{\pi}(\mu) .
\end{equation} 

We now define the set of alternate environments as $K(\mu_L^\pi(\mu), \mu_R^\pi(\mu)) = K_1(\mu_L^\pi(\mu)) \cup K_2(\mu_R^\pi(\mu))$, where
\[
K_1(\mu_L^\pi(\mu)) = \{ \tilde{\nu} : \tilde{\nu} \in \mathbf{S}, \tilde{\mu} < \mu_L^\pi(\mu) \}
\]
and
\[
K_2(\mu_R^\pi(\mu)) = \{ \tilde{\nu} : \tilde{\nu} \in \mathbf{S}, \tilde{\mu} > \mu_R^\pi(\mu) \}.
\]
It is worth noting that to define alternate environments, we also need to choose the cost distribution. We choose the cost distribution to be the same as $\mathcal{C}$, with mean $\overline{c}$, in the alternate environments.

Observe that $\tau_\delta +1$ is a stopping time.  Using the data processing inequality and Wald's lemma (see Lemma 0.1 in \cite{kaufmann2020contributions}), for a given $\delta \in (0,1)$ and any alternate environment $\tilde{\nu}$ with mean $\tilde{\mu}$, we have, using the independence of $\nu$ and $\tilde{\nu}$ with respect to $\mathcal{C}$, that:

	\begin{equation} \nonumber
\mathbb{E}[\tau_\delta +1] \cdot d(\mu, \tilde{\mu})  \geq \sup_{\mathcal{E}_\delta \in \mathcal{F}_{\tau_\delta+1}} \phi(\mathbb{P}_{\nu}(\mathcal{E}_\delta),
	\mathbb{P}_{\tilde{\nu}}(\mathcal{E}_\delta).
  	\end{equation}
It can be re-written as, 

$$\frac{C_\delta}{\log(1/\delta)}\frac{\mathbb{E}[\tau_\delta +1] }{C_\delta}\cdot d(\mu, \tilde{\mu})  \geq \frac{\sup_{\mathcal{E}_\delta \in \mathcal{F}_{\tau_\delta+1}} \phi(\mathbb{P}_{\nu}(\mathcal{E}_\delta),
	\mathbb{P}_{\tilde{\nu}}(\mathcal{E}_\delta))}{\log(1/\delta)}.$$
   Choose  $\mathcal{E_\delta} = \{\tilde{\mu} \notin [\widehat{\mu}^{\pi}_L(\tau_\delta,\delta), \widehat{\mu}^{\pi}_R(\tau_\delta,\delta)]\}$. Observe that, $\mathcal{E}_\delta$ is a measurable event in $\mathcal{F}_{\tau_\delta} \subseteq \mathcal{F}_{\tau_\delta+1}$.

   Using \eqref{appe_10} and the arguments similar to the proof of Theorem 1, we get the desired result.

\hfill $\Box$

\noindent \textbf{Proof of Theorem \ref{thm_kl_cost}.}

First, we show that 
\begin{equation}  \label{appen_12}
     \mathbb{P}_{\nu}(\mu \notin [\widehat{\mu}^{\hat{\pi}_1}_L(\tau_\delta,\delta), \widehat{\mu}^{\hat{\pi}_1}_R(\tau_\delta,\delta)]) \leq \delta.
\end{equation}
 
Using \cite{kaufmann2021mixture}, we get,

\begin{equation} \label{appen_11}
    \mathbb{P}_{\nu} (\exists n \in \mathbb{N}:  \,n d(\hat{\mu}_n,\mu) \geq \beta(n,\delta)) \leq \delta,
\end{equation}
where $\beta(n,\delta)$ is defined in Section 6.
Observe that, 
$$\{ \mu \notin [\widehat{\mu}^{\hat{\pi}_1}_L(\tau_\delta,\delta), \widehat{\mu}^{\hat{\pi}_1}_R(\tau_\delta,\delta)] \} \subseteq \{  \tau_\delta d(\hat{\mu}_{\tau_\delta},\mu) \geq \beta(\tau_\delta,\delta)\}.$$
Further, the following holds as well,

$$\{  \tau_\delta d(\hat{\mu}_{\tau_\delta},\mu) \geq \beta(\tau_\delta,\delta)\} \subseteq \{  \exists n \in \mathbb{N}:  \,n d(\hat{\mu}_n,\mu) \geq \beta(n,\delta)\}.$$

Using \eqref{appen_11}, we get that \eqref{appen_12} holds.

 Observe that from the definition of $\widehat{\mu}^{\hat{\pi}_1}_L(\tau_\delta,\delta)$, $\widehat{\mu}^{\hat{\pi}_1}_R(n,\delta)$ and the strict quasi-convexity of $d(\mu,\cdot)$, we get,

$$d(\hat{\mu}_{\tau_\delta}, \widehat{\mu}^{\hat{\pi}_1}_L(\tau_\delta,\delta)) = d(\hat{\mu}_{\tau_\delta}, \widehat{\mu}^{\hat{\pi}_1}_R(\tau_\delta,\delta)) = \frac{\log(2/\delta)}{\tau_\delta}.$$

It can be re-written as,

$$d(\hat{\mu}_{\tau_\delta}, \widehat{\mu}^{\hat{\pi}_1}_L(\tau_\delta,\delta)) = d(\hat{\mu}_{\tau_\delta}, \widehat{\mu}^{\hat{\pi}_1}_R(\tau_\delta,\delta)) = \frac{\log(2/\delta)}{C_\delta}\frac{C_\delta}{\tau_\delta}.$$

Taking $\delta \to 0$ and using the joint-continuity of $d(\mu,x)$ in $(\mu,x)$ and \eqref{appe_10}, we get that $\pi_1 \in \hat{\Pi}_{\rm CI}^{s}$ and part (b), (c) of this theorem holds. This completes the proof.

\hfill $\Box$

\section{Proofs of results in Section \ref{sec_generalization}.} \label{appen_generalization}

\noindent \textbf{Proof of Theorem \ref{thm_generalization_lower}.}

We prove the result for $\nu \in \mathbf{B}$. A similar proof holds for the case when $\nu \in \mathbf{H}$.

Recall that,
\begin{equation} 
    \tau_\delta = \sup \{ n \in \mathbb{Z}^{+}: \sum_{i = 1}^{n}c_i \leq C_\delta \}. \nonumber
\end{equation}

Again, using elementary renewal theorem for the renewal process, as $C_\delta \to \infty$, we get,

\begin{equation} \label{appe_13}
    \lim_{\delta \to 0}\frac{\mathbb{E}[\tau_\delta]}{C_\delta} = \frac{1}{\overline{c}} \textrm{ almost surely}.
\end{equation}

Consider any $\pi \in \hat{\Pi}_{\rm CI}^{s}$. As the policy $\pi$ is stable, it follows that for a given distribution $\nu \in \mathbf{B}$ with mean $\mu$ and a cost distribution $\mathcal{C}$ with mean $\overline{c}$, we have,
\begin{equation} \nonumber
\lim_{\delta \to 0} \hat{\mu}_L^{\pi}(\tau_\delta,\delta) \overset{p}{\to}  \mu_{L}^{\pi}(\nu) \textrm{ and } \lim_{\delta \to 0} \hat{\mu}_R^{\pi}(\tau_\delta,\delta) \overset{p}{\to}  \mu_{R}^{\pi}(\nu) .
\end{equation} 

We now define the set of alternate environments as $K(\mu_L^\pi(\nu), \mu_R^\pi(\nu)) = K_1(\mu_L^\pi(\nu)) \cup K_2(\mu_R^\pi(\nu))$, where
\[
K_1(\mu_L^\pi(\nu)) = \{ \tilde{\nu} : \tilde{\nu} \in \mathbf{B}, m(\tilde{\nu}) < \mu_L^\pi(\nu) \}
\]
and
\[
K_2(\mu_R^\pi(\nu)) = \{ \tilde{\nu} : \tilde{\nu} \in \mathbf{B}, m(\tilde{\nu}) > \mu_R^\pi(\nu) \}.
\]
It is worth noting that to define alternate environments, we also need to choose the cost distribution. We choose the cost distribution to be the same as $\mathcal{C}$, with mean $\overline{c}$, in the alternate environments.

Observe that $\tau_\delta +1$ is a stopping time. Using the data processing inequality and Wald's lemma, for a given $\delta \in (0,1)$ and any alternate environment $\tilde{\nu}$ with mean $\tilde{\mu}$, we have, using the independence of $\nu$ and $\tilde{\nu}$ with respect to $\mathcal{C}$:

	\begin{equation} \nonumber
\mathbb{E}[\tau_\delta +1] \cdot KL(\nu, \tilde{\nu})  \geq \sup_{\mathcal{E} \in \mathcal{F}_{\tau_\delta+1}} \phi(\mathbb{P}_{\nu}(\mathcal{E}),
	\mathbb{P}_{\tilde{\nu}}(\mathcal{E})).
  	\end{equation}

   First we choose $\tilde{\nu} \in K_1(\mu_L^{\pi}(\nu)) $. Observe that $\inf_{\tilde{\nu} \in K_1(\mu_L^{\pi}(\nu))} KL(\nu,\tilde{\nu}) = KL_{\rm inf}(\nu, \mathbf{B}, \mu_L^{\pi}(\nu))$.
Hence, it can be re-written as, 

$$\frac{C_\delta}{\log(1/\delta)}\frac{\mathbb{E}[\tau_\delta +1] }{C_\delta}\cdot KL_{\rm inf}(\nu, \mathbf{B}, \mu_L^{\pi}(\nu))  \geq \frac{\sup_{\mathcal{E}_\delta \in \mathcal{F}_{\tau_\delta+1}} \phi(\mathbb{P}_{\nu}(\mathcal{E}_\delta),
	\mathbb{P}_{\tilde{\nu}}(\mathcal{E}_\delta))}{\log(1/\delta)}.$$
   Choose  $\mathcal{E}_\delta = \{\tilde{\mu} \notin [\widehat{\mu}^{\pi}_L(\tau_\delta,\delta), \widehat{\mu}^{\pi}_R(\tau_\delta,\delta)]\}$. Observe that, $\mathcal{E}_\delta$ is a measurable event in $\mathcal{F}_{\tau_\delta} \subseteq \mathcal{F}_{\tau_\delta+1}$. 

   Using \eqref{appe_13} and similar to those in the proof of Theorem \ref{thm_spef_lower_bound_width}, we get the desired result.

\hfill $\Box$

\begin{remark}
    Now we give an equivalent version of Theorem \ref{thm_kl_cost} that holds for the policy $\pi_1^{\rm b}$.
\end{remark}

\begin{theorem} \label{thm_kl_cost_appen_b}
  For $\nu \in \mathbf{B}$ and a cost distribution $\mathcal{C}$ with mean $ 0<\overline{c} < \infty$, the policy $\pi_1^{\rm b}$ has following properties:
    
    a) $\pi_1^{\rm b} \in \hat{\Pi}^{s}_{\rm CI}$.
    
    b) If $\lim_{\delta \to 0} \frac{C_\delta} {\log(1/\delta)}\to k$ for $k \in (0,\infty)$, then we have, $   \mu_{R}^{\pi_1^{\rm b}}(\nu) =  \mu_R^*(\nu, k,\overline{c}),$ and $  \mu_{L}^{\pi_1^{\rm b}}(\nu)  =  \mu_{L}^{*}(\nu,k,\overline{c}) $ where, $ \mu_{L}^{*}(\nu,k, \overline{c}) <\mu $ and $ \mu_{R}^{*}(\nu,k,\overline{c}) > \mu$ uniquely solve \eqref{eqn_kl_symmetric_2}.

c) If $\lim_{\delta \to 0} \frac {C_\delta} {\log(1/\delta)}\to 0$, then we have, 
$  \mu_{R}^{\pi_1^{\rm b} }(\nu) - \mu_{L}^{\pi_1^{\rm b} }(\nu)  = 0$.
\end{theorem}

\noindent \textbf{Proof of Theorem \ref{thm_kl_cost_appen_b}.}

Recall that in the proof of Theorem \ref{thm_kl_cost}, we need two properties of $d(\mu,x)$ function. First is that, $d(\mu,x)$ is a strictly quasi-convex function in $x$. The second one is the joint continuity of $d(\mu,x)$ function in $(\mu,x)$. Furthermore, we require an equivalent concentration bound in this setting, as in \eqref{CI_concen}. At last, we also needed the convergence of $\hat{\mu}_{N_\delta}$ to $\mu$.

 Observe that from the definitions of $\mu^{\pi_1^{\rm b}}_{L} (n,\delta)$, $\mu^{\pi_1^{\rm b}}_{R} (n,\delta)$ and the strict quasi-convexity of $KL_{\rm inf}(\hat{\nu}_n, \mathbf{B}, \cdot)$, we get,

$$ KL_{\rm inf}(\hat{\nu}_{N_\delta}, \mathbf{B}, \mu^{\pi_1^{\rm b}}_{L} (n,\delta)) = KL_{\rm inf}(\hat{\nu}_{N_\delta}, \mathbf{B}, \mu^{\pi_1^{\rm b}}_{R} (n,\delta)) = \frac{ \beta(N_\delta,\delta)}{N_\delta}.$$ 

Using results in Appendix F below, we get that, $KL_{\rm inf}(\nu, \mathbf{B},x)$ is a strictly convex function in $x \in (0,1)$ and $KL_{\rm inf}(\nu, \mathbf{B},x)$ is a jointly continuous function in $(\nu,x)$ for $\nu \in \mathbf{B}$ and $x \in (0,1)$. We also know that, empirical distribution $\hat{\nu}_{N_\delta}$ weakly converges to $\nu$ and $\hat{\nu}_{N_\delta} \in \mathbf{B}$. Further, using \eqref{CS_conce_2}, we know that we get a valid CI in both fixed sample size ($N$) and fixed cost budget ($C_\delta$) setting. 

Now, using the arguments given in the proof of Theorem \ref{thm_kl_cost}, we get the desired result. It is worth noting that one can prove an equivalent version of  Theorem \ref{thm_kl} as well.

\hfill $\Box$

\subsection{Asymptotic optimality of KLinf-based construction of CI  for $\nu \in \mathbf{H}$} \label{appen_generalization_1}

We first define $\mathbf{H}$ formally. For a given $\varepsilon>0$ and $\Gamma>0$. Let
\[
\mathbf{H} := \Big\{\kappa \in \mathbf{P}(\mathbb{R}) : \ \mathbb{E}_{X \sim \kappa}\!\big[\,|X|^{1+\varepsilon}\big] \le \Gamma \Big\}.
\]

Similar to the case of $\nu \in \mathbf{B}$, we utilize the CI constructed in Section 3.3.4 of \cite{agrawal2022bandits}.
Let $\hat{\nu}_n$ denotes the empirical distribution after $n$ samples and $\beta(n,\delta) = 1+  \log\left(\frac{2(1+n)^2}{\delta}\right)$. We denote this method as $\pi_1^{\rm h}$. Let, $\mu^{\pi_1^{\rm h}}_{R} (n,\delta) \triangleq \max \{q > m(\hat{\nu}_n) : 
\quad n \, KL_{\rm inf}(\hat{\nu}_n, \mathbf{H}, q) \leq \beta(n,\delta) \}, $ and $\mu^{\pi_1^{\rm h}}_{L} (n,\delta) \triangleq \min \{q < m(\hat{\nu}_n) : 
\quad n \, KL_{\rm inf}(\hat{\nu}_n, \mathbf{H}, q) \leq \beta(n,\delta) \}.$ 
It follows that the reported CI is $[\mu^{\pi_1^{\rm h}}_{L} (n,\delta),\mu^{\pi_1^{\rm h}}_{R} (n,\delta)]$. 

The underlying concentration bound for the above CI construction (see Section 3.3.4 of \cite{agrawal2022bandits}) is 

\begin{equation} \label{new_ai}
      \mathbb{P}_{\nu} ( \exists \, n \in \mathbb{N}: \,n \, KL_{\rm inf}(\hat{\nu}_n, \mathbf{H},\mu) \geq \beta(n,\delta)) \leq \delta.
\end{equation}

\begin{remark}
    Now we give an equivalent version of Theorem \ref{thm_kl_cost} that holds for the policy $\pi_1^{\rm h}$.
\end{remark}

\begin{theorem} \label{thm_kl_cost_appen}
  For $\nu \in \mathbf{H}$ with $\mathbb{E}_{X \sim \nu}[|X|^{1+\varepsilon}] < \Gamma $, the policy $\pi_1^{\rm h}$ has following properties:
    
    a) $\pi_1^{\rm h} \in \hat{\Pi}^{s}_{\rm CI}$.
    
    b) If $\lim_{\delta \to 0} \frac{C_\delta} {\log(1/\delta)}\to k$ for $k \in (0,\infty)$, then we have, $   \mu_{R}^{\pi_1^{\rm h}}(\nu) =  \mu_R^*(\nu, k,\overline{c}),$ and $  \mu_{L}^{\pi_1^{\rm h}}(\nu)  =  \mu_{L}^{*}(\nu,k,\overline{c}) $ where, $ \mu_{L}^{*}(\nu,k, \overline{c}) <\mu $ and $ \mu_{R}^{*}(\nu,k,\overline{c}) > \mu$ uniquely solve \eqref{eqn_kl_symmetric_2}.

c) If $\lim_{\delta \to 0} \frac {C_\delta} {\log(1/\delta)}\to 0$, then we have, 
$  \mu_{R}^{\pi_1^{\rm h} }(\nu) - \mu_{L}^{\pi_1^{\rm h} }(\nu)  = 0$.
\end{theorem}

\noindent \textbf{Proof of Theorem \ref{thm_kl_cost_appen}.}

Using results in Appendix G, we get that $KL_{\rm inf}(\nu, \mathbf{H},x)$ is a strictly convex function in $x \in (0,1)$ and $KL_{\rm inf}(\nu, \mathbf{H},x)$ is a jointly continuously function in $(\nu,x)$ for $\nu \in \mathbf{H}$ and $x \in (0,1)$. We also know that, empirical distribution $\hat{\nu}_n$ weakly converges to $\nu$. It is worth noting that the empirical distribution, denoted by $\hat{\nu}_n$, may not always belong to $\mathbf{H}$. However, since the theorem assumes that $\mathbb{E}_{X \sim \nu}[|X|^{1+\varepsilon}] < \Gamma$, along each sample path, $\hat{\nu}_n \in  \mathbf{H}$ will eventually hold.

To see this, observe that
\[
\mathbb{E}_{X \sim \hat{\nu}_n}[|X|^{1+\varepsilon}] = \frac{1}{n}\sum_{i=1}^{n}|X_i|^{1+\varepsilon}.
\]
By the strong law of large numbers, we have that, for almost every sample path,
\[
\frac{1}{n}\sum_{i=1}^{n}|X_i|^{1+\varepsilon} \longrightarrow \mathbb{E}_{X \sim \nu}[|X|^{1+\varepsilon}] < \Gamma.
\]
Hence, for each sample path, there exists an $N$ such that for all $n \geq N$,
\[
\mathbb{E}_{X \sim \hat{\nu}_n}[|X|^{1+\varepsilon}] \leq \Gamma.
\]
Therefore, along every sample path, the empirical distribution $\hat{\nu}_n$ eventually lies in $\mathbf{H}$.
It is worth noting that the above argument is needed because we only have continuity of $KL_{\rm inf}(\nu,\mathbf{H},x)$ in $\nu$ for $\nu \in \mathbf{H}$. Further, using \eqref{new_ai}, we know that we get a valid CI in both fixed sample size ($N$) and fixed cost budget ($C_\delta$) setting.

Now using the arguments given in the proof of Theorem \ref{thm_kl_cost}, we get the desired result. It is worth noting that one can prove an equivalent version of  Theorem \ref{thm_kl} as well.

\section{Properties of $KL_{\rm inf}(\nu, \mathbf{B},x)$ }
\label{appen_bounded}
For $\nu \in \mathbf{B}$, we define
\begin{align*}
KL_{\rm inf}^{U} (\nu,\mathbf{B},x) &=  \inf_{\kappa \in \mathbf{B} : m(\kappa)\geq x} KL (\nu,\kappa), \quad \text{and}\\
KL_{\rm inf}^{L} (\nu,\mathbf{B},x) &=  \inf_{\kappa \in \mathbf{B} : m(\kappa)\leq x} KL (\nu,\kappa).
\end{align*}

Recall that we had defined $KL_{\rm inf} (\nu,\mathbf{B},x)$ in Section 7. 
Since $KL_{\rm inf}^{U} (\nu,\mathbf{B},x) = 0$ if $ m(\nu) \geq x$ and $KL_{\rm inf}^{L} (\nu,\mathbf{B},x) = 0$ if $ m(\nu) \leq x$, it follows that
\begin{equation} \label{main_sup_84}
     KL_{\rm inf} (\nu,\mathbf{B},x) =  \max\{KL_{\rm inf}^{L} (\nu,\mathbf{B},x), KL_{\rm inf}^{U} (\nu,\mathbf{B},x)\}.
\end{equation}

We now state some properties of $KL_{\rm inf}^{U} (\nu,\mathbf{B},x)$ and $KL_{\rm inf}^{L} (\nu,\mathbf{B},x)$.

\bigskip
\noindent \textbf{Dual Representation of $KL_{\rm inf}^{U} (\nu,\mathbf{B},x)$ and $KL_{\rm inf}^{L} (\nu,\mathbf{B},x)$:} 
In the literature, it is well established that dual representations of $KL_{\rm inf}^{U} (\nu,\mathbf{B},x)$ and $KL_{\rm inf}^{L} (\nu,\mathbf{B},x)$ are much more tractable. 
We now rewrite Theorem~3 of \citet{jourdan2022top}. 
For $(\lambda, \nu, x) \in  \mathbb{R} \times \mathbf{B} \times [0, 1]$, define
\[
H^{+}(\lambda, \nu, x) = \mathbb{E}_{\nu} \big[\log(1-\lambda(X-x))\big],
\]
where we define $\log(z) = -\infty$ for $z \leq 0$. 
Similarly, define
\[
H^{-}(\lambda, \nu, x) = \mathbb{E}_{\nu} \big[\log(1+\lambda(X-x))\big].
\]

\begin{theorem} \label{thm_kl_inf_bounded}
 For all $\nu \in \mathbf{B}$ and $x \in (0,1)$, we have
 $$KL_{\rm inf}^{U} (\nu,\mathbf{B},x)=  \sup_{\lambda \in [0,1/(1-x)]} H^{+}(\lambda, \nu, x), \textrm{ and}$$
  $$KL_{\rm inf}^{L} (\nu,\mathbf{B},x)=  \sup_{\lambda \in [0,1/x]} H^{-}(\lambda, \nu, x).$$
\end{theorem}
It is well-known in the literature that in the above dual representation, it is a univariate convex optimization problem and can be computed efficiently by iterative methods such as Newton's method (see \cite{honda2010asymptotically}).
 
Now we re-write some properties of $KL_{\rm inf}^{U} (\nu,\mathbf{B},x)$ and $KL_{\rm inf}^{L} (\nu,\mathbf{B},x)$ functions which are proven in \citet{honda2010asymptotically} and \citet{jourdan2022top}.
It is worth noting that for continuity of $KL_{\rm inf}^{U} (\nu,\mathbf{B},x)$ and $KL_{\rm inf}^{L} (\nu,\mathbf{B},x)$ in $\mathbf{B}$, we endow it with the topology of weak convergence.
\noindent \subsection{Properties of $KL_{\rm inf}^{U} (\nu,\mathbf{B},x)$ and $KL_{\rm inf}^{L} (\nu,\mathbf{B},x)$ :} \label{KL_inf_properties}
\begin{enumerate} 
    \item  The function $KL_{\rm inf}^{U} (\nu,\mathbf{B},x)$ (resp. $KL_{\rm inf}^{L} (\nu,\mathbf{B},x)$) is continuous on $\mathbf{B} \times [0, 1)$ (resp. $\mathbf{B} \times (0, 1]$). 
    
    \item The function $x \mapsto KL_{\rm inf}^{U}(\nu,\mathbf{B},x)$ is strictly convex on $(m(\nu), 1]$. Further, the function $x \mapsto KL_{\rm inf}^{L}(\nu,\mathbf{B},x)$ is strictly convex on $[0, m(\nu))$. 
\end{enumerate}
 The continuity in the above properties is in this topology on $\mathbf{B}$.

\section{Properties of $KL_{\rm inf}(\nu, \mathbf{H},x)$ }
\label{appen_heavy}
Let \[M := \big[-\Gamma^{\frac{1}{1+\varepsilon}}, \ \Gamma^{\frac{1}{1+\varepsilon}}\big], \ \ M^\circ=\mathrm{int}(M).
\]
For $\nu \in \mathbf{H}$, we define
\begin{align*}
KL_{\rm inf}^{U} (\nu,\mathbf{H},x) &=  \inf_{\kappa \in \mathbf{H} : m(\kappa)\geq x} KL (\nu,\kappa), \quad \text{and}\\
KL_{\rm inf}^{L} (\nu,\mathbf{H},x) &=  \inf_{\kappa \in \mathbf{H} : m(\kappa)\leq x} KL (\nu,\kappa).
\end{align*}

Recall that we had defined $KL_{\rm inf} (\nu,\mathbf{H},x)$ in Section~7. 
Since $KL_{\rm inf}^{U} (\nu,\mathbf{H},x) = 0$ if $ m(\nu) \geq x$ and $KL_{\rm inf}^{L} (\nu,\mathbf{H},x) = 0$ if $ m(\nu) \leq x$, it follows that
\begin{equation} \label{main_sup_H}
     KL_{\rm inf} (\nu,\mathbf{H},x) =  \max\{KL_{\rm inf}^{L} (\nu,\mathbf{H},x), KL_{\rm inf}^{U} (\nu,\mathbf{H},x)\}.
\end{equation}

\noindent\textbf{Dual functions.} For $x\in M^\circ$, $\lambda=(\lambda_1,\lambda_2)\in \mathbb{R}^2$, $\gamma=(\gamma_1,\gamma_2)\in \mathbb{R}^2$, and $X\in\mathbb{R}$, define
\[
g^{U}(X,\lambda,x) \ :=\ 1 \;-\; \lambda_1 (X-x) \;-\; \lambda_2\big(\Gamma - |X|^{1+\varepsilon}\big),
\]
\[
g^{L}(X,\gamma,x) \ :=\ 1 \;+\; \gamma_1 (X-x) \;-\; \gamma_2\big(\Gamma - |X|^{1+\varepsilon}\big).
\]

\noindent\textbf{Feasible dual sets.}
\[
S^U_\gamma(x) \;:=\; \Big\{ \lambda_1\!\ge\!0,\ \lambda_2\!\ge\!0:\ 
1+\lambda_1 x-\lambda_2 \gamma
\;-\;
\frac{\varepsilon\, \lambda_1^{1+\frac{1}{\varepsilon}}}{(1+\varepsilon)^{1+\frac{1}{\varepsilon}}\, \lambda_2^{\frac{1}{\varepsilon}}}
\;\ge\; 0 \Big\},
\]
\[
S^L_\gamma(x) \;:=\; \Big\{ \gamma_1\!\ge\!0,\ \gamma_2\!\ge\!0:\ 
1-\gamma_1 x-\gamma_2 \gamma
\;-\;
\frac{\varepsilon\, \gamma_1^{1+\frac{1}{\varepsilon}}}{(1+\varepsilon)^{1+\frac{1}{\varepsilon}}\, \gamma_2^{\frac{1}{\varepsilon}}}
\;\ge\; 0 \Big\}.
\]
We now re-write the Theorem 4.5 of \cite{agrawal2022bandits}.
\begin{theorem}\label{thm:dual_Lgamma}
Let $\nu\in\mathbf{H}$ and $x\in M^\circ$. Then
\[
KL_{\rm inf}^{U}(\nu,\mathbf{H},x)
\;=\;
\max_{\lambda \in S^U_\gamma(x)}
\ \mathbb{E}_\nu\!\left[\log\!\big(g^U(X,\lambda,x)\big)\right],
\]
\[
KL_{\rm inf}^{L}(\nu,\mathbf{H},x)
\;=\;
\max_{\gamma \in S^L_\gamma(x)}
\ \mathbb{E}_\nu\!\left[\log\!\big(g^L(X,\gamma,x)\big)\right].
\]

\end{theorem}

The above dual representation is a bivariate convex optimization problem and can be computed efficiently by iterative methods.
We now restate the relevant parts of Lemma 4.9, Lemma 4.10 and Lemma 4.11 of \cite{agrawal2022bandits}.

\begin{lemma}
 For $\eta \in \mathbf{H}$ and $x \in M^{o}$ such that $x > m(\eta)$, $KL_{\mathrm{inf}}^{U}(\eta, \mathbf{H}, x)$ is a strictly convex function of $x$. 
\end{lemma}

\begin{lemma}
 For $\eta \in \mathbf{H}$ and $x \in M^{o}$ such that $x < m(\eta)$, $KL_{\mathrm{inf}}^{L}(\eta, \mathbf{H}, x)$ is a strictly convex function of $x$. 
\end{lemma}

\begin{lemma}
When restricted to $\mathcal{H} \times M^{o}$, $KL_{\mathrm{inf}}^{U}$ and $KL_{\mathrm{inf}}^{L}$ are jointly continuous in their arguments.
\end{lemma}

 The continuity in the above Lemma is in the topology of weak convergence on $\mathbf{H}$.

\section{Numerical study for the case when $\nu \in \mathbf{H}$} \label{appen_num}
We now simulate our policy $\pi_1^{h}$. The experiment is conducted on Pareto distributions with the scale parameter $x_m = 1$ and the shape parameter $\alpha = 3$. For the description of $\mathbf{H}$, we have chosen $\epsilon = 1$ and $\Gamma = 4$.
We run the experiment in the set-up where we have costly samples with cost budget $C \in \{500, 1000, 2000, 3000\} $.
 We assume the cost distribution to be of uniform distributon in $[0,2]$, i.e., $\mathcal{C}= Unif[0,2]$ and $\delta$ is set to be $5 \%$.

We generated 1000 i.i.d. datasets of size \( C \in \{2000, 3000\} \), computed CI widths, and report the average width (max 95$\%$ CI width: $0.01$) in Table \ref{tab:ci_widths}.

\begin{table}[t]
\centering
\caption{Average CI widths over 1000 i.i.d. datasets.}
\label{tab:ci_widths}
\begin{tabular}{lc}
\hline
\textbf{Cost Budget \(C\)} & \textbf{Avg.\ CI width}  \\
\hline
500 & 0.492 \\
1000 & 0.355   \\ 
2000 & 0.255 \\
3000 & 0.199 \\
\hline
\end{tabular}
\end{table}

\section{One-sided CI setting} \label{one_sided_CI}

Let $\Pi^s_{\rm CI_L}$ denote the collection of stable policies for constructing one-sided CIs of the form $[\widehat{\mu}^{\pi}_L(N, \delta), \overline{\mu})$ for the mean after observing  $X_1, X_2, \ldots, X_N$ for any $\delta \in (0,1)$. For any policy $\pi \in \Pi^s_{\rm CI_L}$ and a given distribution $\nu $ with mean $\mu$, the following condition must hold for any $\delta \in (0,1)$:
$
\forall n \in \mathbb{N} : \mathbb{P}_{\nu}(\mu \in [\widehat{\mu}^{\pi}_L(n, \delta), \overline{\mu})) \geq 1-\delta, \textrm{and }
$ $
\lim_{\delta \to 0} \widehat{\mu}^{\pi}_L(N_\delta, \delta) \overset{p}{\to} \mu_{L}^{\pi}(\nu),
$
where $\mu_{L}^{\pi}(\nu) \leq \mu$ is a constant. Here, we define $\mu - \mu_{L}^{\pi}(\nu) $ as the ``half-width'' of the one-sided CI of the mean. Now we state the result that characterizes the three learning regimes for this setting.
\begin{corollary} \label{cor_one_sided_CI}
For a given $\nu \in \mathbf{S}$ with mean $\mu$, and any $\pi \in \Pi^s_{\rm CI_L}$, the following holds:
\begin{enumerate}
    \item \textbf{No learning regime:} If $\lim_{\delta \to 0} \frac{N_\delta} {\log(1/\delta)} \to 0$, then $
    \lim_{\delta \to 0} [\mu - \widehat{\mu}^{\pi}_L(N_\delta, \delta) ]\overset{p}{\to} \mu- \underline{\mu}.$

    \item \textbf{Sufficient learning regime:} If $\lim_{\delta \to 0} \frac{N_\delta} {\log(1/\delta)} \to k$ for $k \in (0, \infty)$, then
 $ \mu- \mu_{L}^{\pi}(\mu)  \geq \mu-  \mu_{L}^*(\mu, k), $
    where  $\mu_{L}^*(\mu, k) < \mu$ uniquely solves $  d(\mu, \mu_{L}^*(\mu, k))= \frac{1}{k}.$
\end{enumerate}
\end{corollary}

The above results can similarly be applied to one-sided CIs where the goal is to ensure the upper bound of the interval is below a threshold by replacing all occurrences of subscripts $L$ with $R$ in the definitions and results.


\noindent \textbf{Proof of Corollary \ref{cor_one_sided_CI}}.

Proof of this corollary follows similar to the proof of Theorem \ref{thm_spef_lower_bound_width}.

\begin{remark}
    It is worth noting that $\pi_1$ can be used to construct an asymptotically optimal one-sided CI. Observe that the one-sided CI using $\pi_1$ is given by $[\hat{\mu}_{L}^{\pi_1}(N,\delta), \overline{\mu})$, where $\hat{\mu}_{L}^{\pi_1}(N,\delta)$ is defined in \eqref{policy_pi_1}. 
\end{remark}

\begin{remark}
   The results for one-sided CI can be extended for the non-parametric and costly sample cases as well trivially.
\end{remark}

\section{Non-asymptotic analysis of CI Width for our policies}
\label{appen_finite_confidence}

We provide a non-asymptotic analysis of CI width for our policies for both the canonical exponential-family case and the non-parametric (bounded) distributional case.
\paragraph{Exponential family case.}
Fix $\nu \in \mathbf{S}$ and let $\mu$ be its mean. Let the variance function be
$\sigma^2(\mu)$.

Recall, under policy $\pi_1$, we have
\[
    d\!\left(\hat{\mu}_n,\, \hat{\mu}_R^{\pi_1}(n,\delta)\right)
    = d\!\left(\hat{\mu}_n,\, \hat{\mu}_L^{\pi_1}(n,\delta)\right)
    = \frac{\log(2/\delta)}{n}.
\]
Using a Taylor-series expansion of $d(\mu, x)$ around $x = \mu$, we get
\[
    \frac{\left(\hat{\mu}_R^{\pi_1}(n,\delta) - \hat{\mu}_n\right)^2}{2\,\sigma^2(c_{n,R})}
    = \frac{\left(\hat{\mu}_n - \hat{\mu}_L^{\pi_1}(n,\delta)\right)^2}{2\,\sigma^2(c_{n,L})}
    = \frac{\log(2/\delta)}{n},
\]
for some $c_{n,R} \in \bigl[\hat{\mu}_n,\, \hat{\mu}_R^{\pi_1}(n,\delta)\bigr]$ and
$c_{n,L} \in \bigl[\hat{\mu}_L^{\pi_1}(n,\delta),\, \hat{\mu}_n\bigr]$.

Hence, we obtain the width
\[
    \hat{\mu}_R^{\pi_1}(n,\delta) - \hat{\mu}_L^{\pi_1}(n,\delta)
    = \sqrt{\frac{2\log(2/\delta)}{n}}\,
      \Bigl(\sigma(c_{n,R}) + \sigma(c_{n,L})\Bigr).
\]
Assume the mean space contains a compact interval $I$ such that
$\hat{\mu}_L^{\pi_1}(n,\delta),\, \hat{\mu}_n,\, \hat{\mu}_R^{\pi_1}(n,\delta) \in I$.
Define $\sigma_{\max} = \sup_{\mu \in I} \sigma(\mu) < \infty$. Then the width is bounded
by
\[
    \hat{\mu}_R^{\pi_1}(n,\delta) - \hat{\mu}_L^{\pi_1}(n,\delta)
    \;\leq\; 2\,\sigma_{\max} \sqrt{\frac{2\log(2/\delta)}{n}}.
\]

\paragraph{Non-parametric (bounded) case.}
Now, we perform a similar analysis for fixed $\nu \in \mathbf{B}$. Let $m(\nu)$ denote
the mean of $\nu$. Recall, under policy $\pi^{b}_1$, we have
\[
    n\, KL_{\rm inf}\!\left(\hat{\nu}_n,\, \mathbf{B},\, \hat{\mu}_R^{\pi_1^b}(n,\delta)\right)
    = n\, KL_{\rm inf}\!\left(\hat{\nu}_n,\, \mathbf{B},\, \hat{\mu}_L^{\pi_1^b}(n,\delta)\right)
    = \beta(n,\delta).
\]
Using Proposition~1 in \cite{orabona2023tight} and Pinsker's inequality, we have

\[
    \frac{1}{n}\,\beta(n,\delta)
    = KL_{\rm inf}\!\left(\hat{\nu}_n,\, \mathbf{B},\, \hat{\mu}_R^{\pi_1^b}(n,\delta)\right)
    \;\geq\; 2\,\Bigl(\hat{\mu}_R^{\pi_1^b}(n,\delta) - \hat{\mu}_n\Bigr)^2.
\]
A similar argument holds for $x = \hat{\mu}_L^{\pi_1^b}(n,\delta)$. Hence,
\[
    \hat{\mu}_R^{\pi_1^b}(n,\delta) - \hat{\mu}_L^{\pi_1^b}(n,\delta)
    = \Bigl(\hat{\mu}_R^{\pi_1^b}(n,\delta) - \hat{\mu}_n\Bigr)
      + \Bigl(\hat{\mu}_n - \hat{\mu}_L^{\pi_1^b}(n,\delta)\Bigr)
    \;\leq\; \sqrt{\frac{2\,\beta(n,\delta)}{n}}.
\]
Using $\beta(n,\delta) = 1 + \log\!\left(\frac{2(1+n)}{\delta}\right)$, we obtain the
explicit bound
\[
    \hat{\mu}_R^{\pi_1^b}(n,\delta) - \hat{\mu}_L^{\pi_1^b}(n,\delta)
    \;\leq\; \sqrt{\frac{2\!\left(1 + \log\!\left(\dfrac{2(1+n)}{\delta}\right)\right)}{n}}.
\]

\section{Heuristic connections of Cram\'er--Rao lower bound to our lower bound} \label{appen_cramer}
In this section, we draw a heuristic connection between the Cramér--Rao lower bound and our lower bound in the canonical one-parameter exponential family setting. In a canonical one-parameter
exponential family, i.e., $\mathbf{S}$ parameterized by its mean $\mu$, the Fisher information equals
\[
    I(\mu) = \frac{1}{\sigma^2(\mu)},
\]
where $\sigma^2(\mu)$ denotes the variance of the true distribution $\nu$ with mean $\mu$.
Accordingly, the KL divergence admits the local expansion
\[
    d(\mu, \mu') = \frac{1}{2} I(\mu)(\mu' - \mu)^2 + o\!\left((\mu' - \mu)^2\right).
\]
It follows that
\[
    d(\mu, \mu') = \frac{1}{2\sigma^2(\mu)}(\mu' - \mu)^2 + o\!\left((\mu' - \mu)^2\right).
\]
In the sufficient learning regime, the optimal endpoints satisfy
\[
    d\!\left(\mu,\, \mu_L^*(\mu, k)\right) = d\!\left(\mu,\, \mu_R^*(\mu, k)\right) = \frac{1}{k},
\]
and the quadratic approximation yields
\[
    \left|\mu_R^*(\mu, k) - \mu\right| \approx \left|\mu - \mu_L^*(\mu, k)\right|
    \approx \sqrt{\frac{2\sigma^2(\mu)}{k}},
\]
so that
\[
    \mu_R^*(\mu, k) - \mu_L^*(\mu, k) \approx 2\sqrt{\frac{2\sigma^2(\mu)}{k}}.
\]
A similar expression arises by combining the Cram\'er--Rao lower bound with the classical
CLT. Let $\hat{\mu}_N$ be any regular estimator for which
\[
    \sqrt{N}\,(\hat{\mu}_N - \mu) \xrightarrow{d} \mathcal{N}(0,\, v(\mu)).
\]
By the Cram\'er--Rao lower bound, $v(\mu) \geq 1/I(\mu) = \sigma^2(\mu)$. In the
sufficient learning regime with $N = k\log(1/\delta)$, a CLT-based CI of half-width
$z_{1-\delta/2}\sqrt{v(\mu)/N}$ has limiting width
\[
    2z_{1-\delta/2}\sqrt{\frac{v(\mu)}{N}}
    \;\approx\; 2\sqrt{\frac{2\log(1/\delta)\,v(\mu)}{k\log(1/\delta)}}
    = 2\sqrt{\frac{2\,v(\mu)}{k}}
    \;\geq\; 2\sqrt{\frac{2\sigma^2(\mu)}{k}},
\]
with equality when $v(\mu) = \sigma^2(\mu)$, i.e., when the estimator is efficient. This
heuristically confirms that the KL-based CI attains the Cram\'er--Rao optimal limiting
width in the sufficient learning regime.

\end{document}